\documentclass[twocolumn]{autart}    

\usepackage{ifpdf}
\usepackage{cite}
\usepackage{amsmath}
\usepackage{amssymb}
\usepackage{amsfonts}
\usepackage{url}
\usepackage{fixltx2e}
\usepackage{graphics} 
\usepackage{graphicx}
\usepackage{epsfig} 
\usepackage{stfloats}
\usepackage{epstopdf}	
\usepackage{apalike}
\usepackage[colorlinks=true,linkcolor=blue, citecolor=blue, urlcolor=blue]{hyperref}
\allowdisplaybreaks

\newcommand{\tr}{\mathrm{tr}}
\newcommand{\T}{^{\top}} 
\newcommand{\dom}{\text{dom }}
\newcommand{\ie}{\textit{i.e.}}

\begin{document}

	\begin{frontmatter}
		
		\title{Nonlinear State Estimation for Inertial Navigation Systems With Intermittent Measurements} 
		
		\thanks[footnoteinfo]{This work was supported by the National Sciences and Engineering Research Council of Canada-NSERC-DG 228465-2013-RGPIN. The material in this paper was not presented at any conferences. Corresponding author A. Tayebi, Tel.: (807) 343-8597, fax: (807) 766-7243.}
		
		\author[UWO]{Miaomiao Wang}\ead{mwang448@uwo.ca}  ~~~~  
		\author[UWO,LU]{Abdelhamid Tayebi}\ead{atayebi@lakeheadu.ca}    

		\address[UWO]{Department of Electrical and Computer Engineering, Western University, London, Ontario, Canada, N6A 3K7
		}  
		\address[LU]{Department of Electrical Engineering, Lakehead University, Thunder Bay, Ontario, Canada, P7B 5E1
		}

		\begin{keyword}                           
			Inertial navigation, nonlinear observers, hybrid dynamical systems, intermittent measurements          
		\end{keyword}                             

		\begin{abstract}                          
			This paper considers the problem of simultaneous estimation of the attitude, position and linear velocity for vehicles navigating in a three-dimensional space. We propose two types of hybrid nonlinear observers using continuous angular velocity and linear acceleration measurements as well as intermittent landmark position measurements. The first type relies on a fixed-gain design approach based on an infinite-dimensional optimization, while the second one relies on a variable-gain design approach based on a continuous-discrete Riccati equation. For each case, we provide two different observers with and without the estimation of the gravity vector. The proposed observers are shown to be exponentially stable with a large domain of attraction. Simulation and experimental results are presented to illustrate the performance of the proposed observers.
		\end{abstract}
		
	\end{frontmatter}
	
	\endNoHyper 
	
	\section{Introduction}
	In the present work, we are interested in the problem of simultaneous estimation of the attitude, position and linear velocity of a rigid body evolving in a three-dimensional space.  This type of estimation, referred to as inertial navigation observer, is crucial for autonomous navigation systems. It is well known that the attitude can be estimated using body-frame observations of some known inertial vectors, \textit{e.g.,} using a star tracker or an inertial measurement unit (IMU), while the position and linear velocity can be obtained, for instance, from a Global Positioning System (GPS). However, it is challenging to design inertial navigation observers for applications in GPS-denied environments (\textit{e.g.,} indoor applications). As an alternative solution to the lack of GPS information, one can use, for instance, bearing measurements from vision system or range measurements from Ultra Wideband systems. In this context, some techniques combining IMU and vision-based landmark position measurements have been proposed in the literature \cite{rehbinder2003pose,mourikis2007multi,mourikis2009vision}. A class of nonlinear invariant pose (attitude and position) observers, designed on the matrix Lie group $SE(3)$ using group velocity (angular and linear velocities) and vision-based landmark position measurements, have been proposed in \cite{vasconcelos2010nonlinear,hua2015gradient,khosravian2015observers} with almost global asymptotic stability guarantees, and in \cite{wang2017globally,wang2019hybrid} with global asymptotic/exponential stability guarantees.		
	
	In practice, obtaining the linear velocity using low-cost sensors in GPS-denied environments is not an easy task. Therefore, it is of great importance to develop estimation algorithms that provide the pose and linear velocity using inertial-vision systems. Vision systems for autonomous navigation have been widely used in robotics applications for many years \cite{kelly2011visual,hesch2013consistency,scaramuzza2011visual}. Most of the existing results in the literature for the pose and linear velocity estimation rely on Kalman-type filters such as the extended Kalman filter (EKF) and unscented Kalman filter (UKF) \cite{mourikis2007multi,mourikis2009vision}. It is well known that these Kalman-type filters, relying on local linearizations, suffer from large computational overhead and lack of strong stability guarantees. Recently, nonlinear geometric observers for inertial navigation systems, using IMU and landmark position measurements, have emerged in the literature. An invariant Extended Kalman Filter (IEKF) has been proposed in  \cite{barrau2017invariant}, and a Riccati-based observer, with gravity estimation, has been proposed in \cite{hua2018riccati}. In contrast to these observers with local stability guarantees,  hybrid nonlinear geometric observers with global exponential stability guarantees have been proposed in \cite{wang2018navigation,wang2020hybrid}.
	
	On the other hand, some interesting results, addressing the estimation problem with intermittent measurements, have been proposed in \cite{ferrante2016state,li2017robust,sferlazza2019time,alonge2019hybrid,berkane2019attitude}. This problem is motivated by the fact that some applications may involve different types of sensors with different bandwidths and communication delays, and as such, irregular sensors sampling may take place. For instance, IMU measurements can be easily considered as continuous compared to visual measurements which often require low sampling rates due to hardware limitations of the vision sensors and the heavy image processing computations. Therefore, the stability is not guaranteed if one tries to implement continuous-time observers in applications involving intermittent measurements combining sensors with different bandwidth characteristics (such as IMU and vision systems), and as such, the observers need to be carefully redesigned.

	In this paper, we propose two types of hybrid nonlinear observers, with fixed and variable gains, for inertial navigation systems, relying on continuous angular velocity and linear acceleration measurements, and intermittent landmark measurements. The main advantage of the variable gain observers is their efficiency in handling measurement noise via a systematic tuning of the gains. For each type, we provide two different observers, with and without the knowledge of the gravity vector, endowed with exponential stability guarantees with a large domain of attraction. The exponential stability results obtained in this paper do not rely on linearizations compared to the recent work in \cite{barrau2017invariant,hamel2018riccati}. In fact, the proposed observers do not have any restrictions on the initial conditions of the position and linear velocity. In contrast to the present work, the hybrid observers proposed in our previous work \cite{wang2020hybrid} are not designed to handle intermittent landmark measurements. Moreover,  our hybrid observers with gravity estimation do not require the knowledge of the gravity vector, which was not considered in \cite{barrau2017invariant,wang2020hybrid}.  Unlike the results of \cite{BerkaneCDC2017,berkane2019attitude}, the estimated attitude from our hybrid observers is continuous, which is desirable in practice, especially when dealing with observer-controller implementations.
	
	The remainder of this paper is organized as follows: Section \ref{sec:preliminary} introduces some preliminary notions that will be used throughout this paper. Section \ref{sec:observerI} is devoted to the design of the nonlinear   observers for inertial navigation systems with fixed-gain design and variable-gain design. Simulation and experimental results are presented in Section \ref{sec:simulation} and Section \ref{sec:experimental}, respectively.

	\section{Preliminary Material}\label{sec:preliminary}
	\subsection{Notations}
	The sets of real, non-negative real, natural numbers and non-zero natural numbers are denoted as $\mathbb{R}$, $\mathbb{R}_{ \geq 0}$, $\mathbb{N}$ and $\mathbb{N}_{>0}$, respectively. We denote by $\mathbb{R}^n$ the $n$-dimensional Euclidean space, and denote by $\mathbb{S}^n$ the set of $(n+1)$-dimensional unit vectors.   The Euclidean norm of a vector $x\in \mathbb{R}^n$ is defined as $\|x\| = \sqrt{x\T x}$, and the Frobenius norm of a matrix $X\in \mathbb{R}^{n\times m}$ is given by $\|X\|_F = \sqrt{\tr(X\T X)}$. A $n$-by-$n$ identity matrix is denoted by $I_n$ and a $n$-by-$m$ zero matrix is denoted by $0_{n\times m}$. For a given matrix $A\in \mathbb{R}^{n\times n}$, we define $\mathcal{E}(A)$ as the set of all unit eigenvectors of $A$. We denote by $\lambda^A_i$ the $i$-th eigenvalue of $A$, and by $\lambda^A_{m}$ and $\lambda^A_{M}$ the minimum and maximum eigenvalue of $A$, respectively.
	
	Let $\{\mathcal{I}\}$ be an inertial frame and $\{\mathcal{B}\}$ be a frame attached to a rigid body. The matrix $R\in SO(3)$ denotes the rotation of frame $\{\mathcal{B}\}$ with respect to frame $\{\mathcal{I}\}$, where $SO(3):=\{R\in \mathbb{R}^{3\times 3}| RR\T = R\T R= I_3, \det{R}=1\}$. We denote by $p\in \mathbb{R}^3$ and $v\in \mathbb{R}^3$, the position and linear velocity of the rigid-body expressed in the inertial frame $\{\mathcal{I}\}$.
	The \textit{Lie algebra} of  $SO(3)$, denoted by $\mathfrak{so}(3)$, is given by
	$\mathfrak{so}(3) = \{\Omega \in \mathbb{R}^{3\times 3}| \Omega = -\Omega\T \}$.
	We denote by $\times$  the vector cross-product on $\mathbb{R}^3$, and define the map $(\cdot)^\times: \mathbb{R}^3 \to \mathfrak{so}(3)$ such that $x\times y = x^\times y, \forall x,y\in \mathbb{R}^3 $.
	For any $R\in SO(3)$, we define $|R|_I\in [0,1]$ as the normalized Euclidean distance on $SO(3)$ with respect to the identity $I_3$, such that $|R|_I^2 = \frac{1}{8}\|I_3-R\|_F^2=\frac{1}{4}\tr(I_3-R)$. Let the map $\mathcal{R}_a: \mathbb{R}\times \mathbb{S}^2$ represent the well-known angle-axis parameterization of the attitude, which is given by $\mathcal{R}_a(\theta,u):=I_3 + \sin\theta u^\times + (1-\cos\theta) (u^\times)^2$ with $\theta \in \mathbb{R}$ denoting the rotation angle and $u\in \mathbb{S}^2$ denoting the rotation axis.
	For any matrix $A=[a_{ij}]_{1\leq i,j\leq 3} \in \mathbb{R}^{3\times 3}$ and any vector $u\in \mathbb{R}^3$, one can verify that $\tr(A u^\times) = -2 u\T \psi(A)$, with $\psi(A)  :=vec(\mathbb{P}_a(A))=\frac{1}{2}[a_{32}-a_{23}, a_{13}-a_{31}, a_{21}-a_{12}]\T$, where $vec (\cdot)$ is the inverse map of $(\cdot)^{\times}$ and $\mathbb{P}_a(A)=(A-A\T)/2$ is the anti-symmetric part of $A$.

	\subsection{Hybrid Systems Framework} 
	Define the \textit{hybrid time domain} as a subset $E \subset \mathbb{R}_{ \geq 0} \times \mathbb{N}$ in the form of
	$ E = \bigcup_{j=0}^{J-1} ([t_j,t_{j+1}] \times \{j\}),	$
	for some finite sequence $0=t_0 \leq t_1 \leq \cdots \leq t_J$, with the last interval possibly in the form $[t_{j}, T)$ with $T$ finite or $T=+\infty$. 
	Given a smooth manifold $\mathcal{M}$ embedded in $\mathbb{R}^n$, define $T\mathcal{M}$ as the tangent space of $\mathcal{M}$. We consider the following hybrid system \cite{goebel2009hybrid}:
	\begin{equation}\mathcal{H}:  ~~~
	\begin{cases}
	\dot{x} ~~\in F(x),& \quad x \in \mathcal{F}   \\
	x^{+} \in G(x),& \quad x \in \mathcal{J}
	\end{cases} \label{eqn:hybrid_system}
	\end{equation}
	where the \textit{flow map} $F: \mathcal{M} \to T\mathcal{M}$ describes the continuous flow of $x$ on the \textit{flow set} $\mathcal{F} \subseteq \mathcal{M}$; the \textit{jump map} $G: \mathcal{M}\rightrightarrows  \mathcal{M}$ (a set-valued mapping from $\mathcal{M}$ to $\mathcal{M}$) describes the discrete flow of $x$ on the \textit{jump set} $\mathcal{J} \subseteq \mathcal{M}$. A hybrid arc is a function $x: \dom{x} \to \mathcal{M}$, where $\dom{x}$ is a hybrid time domain and, for each fixed $j\in \mathbb{N}_{>0}$, $t \mapsto x(t,j)$ is a locally absolutely continuous function on the interval $I_j = \{t:(t,j) \in \dom{x}\}$. Note that $x^+$  denotes the value of $x$ after a jump, namely, $x^+=x(t,j+1)$ with $x(t,j)$ denoting the value of $x$ before the jump. For more details on dynamic hybrid systems, we refer the reader to \cite{goebel2009hybrid,goebel2012hybrid} and references therein.  
	Moreover, we consider the following notion of exponential stability of closed sets for a general hybrid system: a closed set $\mathcal{A}$ is said to be (locally) exponentially stable for the hybrid system $\mathcal{H}$ if there exist strictly positive scalars $\kappa, \lambda$ and $\mu$ such that, for any $|\phi(0,0)|_{\mathcal{A}} < \mu$, every maximal solution to $\mathcal{H}$ is complete  and satisfies 
	$|\phi(t,j)|_{\mathcal{A}} \leq \kappa \exp(-\lambda(t+j))|\phi(0,0)|_{\mathcal{A}}$ for all $(t,j)\in \dom \phi$ with $|\phi|_{\mathcal{A}}$ denoting the distance of $\phi$ to $\mathcal{A}$ \cite{teel2013lyapunov}.

	\subsection{Kinematics and Measurements}
	Consider the following kinematics of a rigid body navigating in a three-dimensional space:
	\begin{align}
	\dot{R} & = R \omega^\times \label{eqn:R}\\
	\dot{p} & = v \label{eqn:p}\\
	\dot{v} & = g  + Ra \label{eqn:v}
	\end{align}
	where $g \in \mathbb{R}^3$ denotes the gravity vector with $c_g = \|g \|$ being the gravity constant,  $\omega\in \mathbb{R}^3$ denotes the angular velocity expressed in the body-frame, and $a\in \mathbb{R}^3$ is the apparent acceleration capturing all non-gravitational forces applied to the rigid body expressed in the body-frame. We assume that the measurements of $\omega$ and $a$ are continuously available.
	
	Consider a family of $N$ landmarks with $p_i \in \mathbb{R}^3$ being the position of the $i$-th landmark expressed in the inertial frame $\{\mathcal{I}\}$. The landmark measurements expressed in the body frame $\{\mathcal{B}\}$ are denoted as
	\begin{equation}
	y_i := R\T(p_i-p), \quad i=1,2,\cdots,N. \label{eqn:output_y}
	\end{equation}
	Note that the landmark measurements can be directly constructed, for instance, from a stereo vision system.
	

	\begin{assum}\label{assum:1}
		The landmark measurements are available at some instants of time $t_j, j\in \mathbb{N}_{>0}$, and there exist constants $0<T_m\leq T_M<\infty$ such that $0\leq t_1\leq T_M$ and $T_m\leq t_{j+1}-t_{j} \leq T_M$ for all $j\in \mathbb{N}_{>0}$.	
	\end{assum}
	This assumption guarantees that the time between two consecutive measurements is lower and upper bounded. The lower bound $T_m$ is required to be strictly positive to avoid Zeno behaviors. Note that in the case where $T_m=T_M$, we have a regular periodic sampling.
	
\begin{figure}[!thpb]
		\centering
		\includegraphics[width=0.9\linewidth]{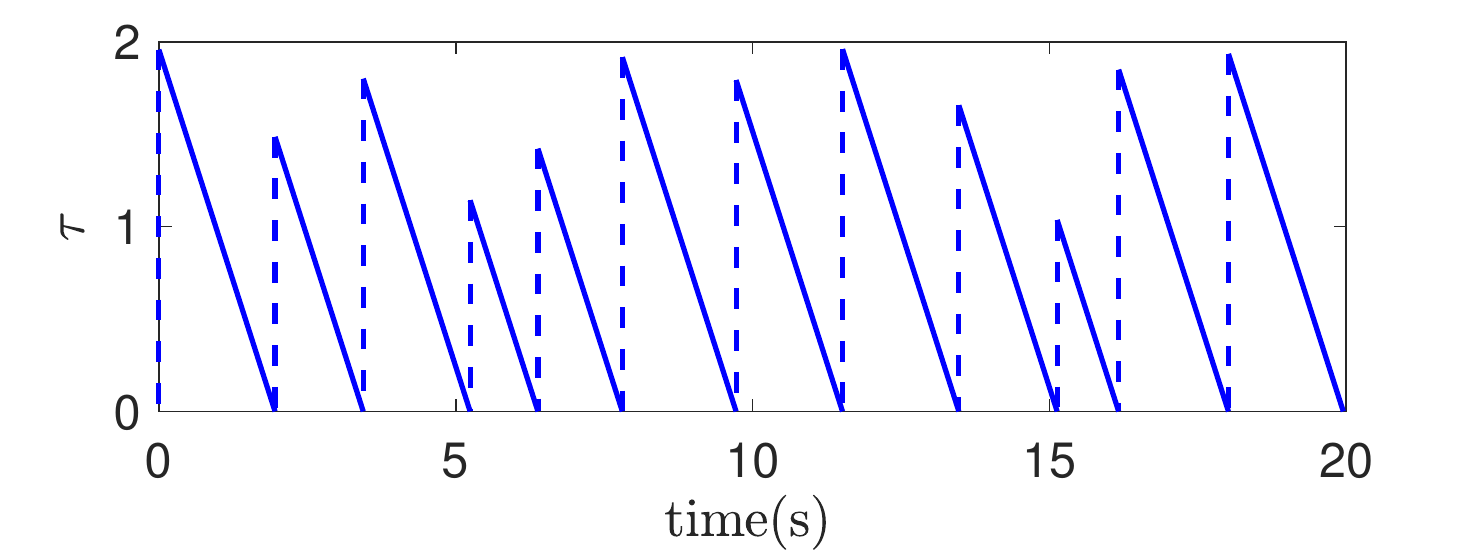}
		\caption{An example of the solution of the timer $\tau$ with $T_m = 1$ and $T_M = 2$.}
		\label{fig:timertau}
	\end{figure}

	To keep track of these time-driven sampling events, a virtual  timer $\tau$, motivated by  \cite{carnevale2007lyapunov,ferrante2016state}, is considered with the following hybrid dynamics:
	\begin{equation}
	\begin{cases}
	\dot{\tau} ~~= -1 & \tau\in [0,T_M] \\
	\tau^+ \in [T_m, T_M] & \tau\in \{0\}
	\end{cases} \label{eqn:tau}
	\end{equation}
	with $\tau(0)\in [0,T_M]$. This virtual state $\tau$ decreases to zero continuously, and upon reaching zero (\ie,  at the arrival of the landmark measurements) it is reset to a value between $T_m$ and $T_M$.
	An example of the solution of the timer $\tau$ is shown in Fig. \ref{fig:timertau}. 
	With this additional state $\tau$ the time-driven sampling events can be described as state-driven events. Note that a different increasing timer has also been used in \cite{carnevale2007lyapunov,sferlazza2019time,berkane2019attitude}. The decreasing timer is purposefully chosen here as it suits our stability proofs that will be given later.
	
	\begin{assum}\label{assum:2}
		There exist at least three non-collinear landmarks among the $N \geq 3$ measurable landmarks at each instant of time $t_j, j\in \mathbb{N}_{>0}$.
	\end{assum}
	Assumption \ref{assum:2} is commonly used in the problem of pose estimation   \cite{vasconcelos2010nonlinear,hua2015gradient,khosravian2015observers,wang2019hybrid} and state estimation for inertial navigation \cite{barrau2017invariant,wang2020hybrid}. Consider a set of scalars $k_i>0, i=1,2,\cdots N$ such that $k_c:= \sum_{i=1}^N k_i=1$. Define $p_c := \sum_{i=1}^{N} k_i p_i$ as the weighted center of landmarks in the inertial frame.  Given three non-collinear landmarks and $k_i>0, \forall i=1,2,\cdots N$, it is always possible to guarantee  that the matrix $M: = \sum_{i=1}^N k_i (p_i-p_c)(p_i-p_c)\T$ is  positive semi-definite with no more than one zero eigenvalue under Assumption \ref{assum:2}.

	Our objective is to design a hybrid estimation scheme for pose and linear velocity using the above mentioned available measurements under Assumption \ref{assum:1} and Assumption \ref{assum:2}.

	\section{Hybrid Observers Design} \label{sec:observerI}
	\subsection{Fixed-gain design}
	\subsubsection*{Known gravity case (without gravity estimation):}
	Let $\hat{R}\in SO(3)$ denote the estimate of the attitude $R$, $\hat{p} \in \mathbb{R}^3$ denote the estimate of the position $p$, and $\hat{v} \in \mathbb{R}^3$ denote the estimate of the linear velocity $v$. We will make use of an auxiliary variable $\eta\in \mathbb{R}^3$ with hybrid dynamics designed to model the intermittent measurements for attitude estimation. This auxiliary variable remains constant between two consecutive landmark measurements  (\ie, $t\in(t_{j-1},t_{j}), \forall j\in \mathbb{N}_{>0}$) and updates upon the arrival of the landmark measurements (\ie, $t=t_j, j\in \mathbb{N}_{>0}$). The estimated attitude $\hat{R}$ is obtained through a continuous integration of the attitude kinematics using the angular velocity $\omega$ and the auxiliary variable $\eta$. The position and velocity are obtained via a hybrid observer consisting of a continuous integration of the translational dynamics using $a$, $g$ and $\eta$ between two consecutive landmark measurements, and a discrete update upon the arrival of the landmark measurements. Our proposed hybrid observer is given as follows:
	\begin{align}
	\underbrace{
		\begin{array}{ll}
		\dot{\hat{R}} &= \hat{R}(\omega + \hat{R}\T \eta )^\times   \\
		\dot{\eta}  &= 0_{3\times 1}  \\
		\dot{\hat{p}} &= \eta^\times (\hat{p}-p_c) + \hat{v}   \\
		\dot{\hat{v}} &= \eta^\times  \hat{v} + g  + \hat{R}a
		\end{array}~}_{\tau \in [0, T_M] }
	\quad
	\underbrace{
		\begin{array}{ll}
		\hat{R}^+ &= \hat{R}  \\
		\eta^+    & = k_R \sigma_R   \\
		\hat{p}^+ &=  \hat{p}  + k_p y     \\
		\hat{v}^+ &=  \hat{v}  + k_v  y
		\end{array}~  }_{\tau \in \{0\}}
	\label{eqn:observer1}
	\end{align}
	where $\hat{R}(0)\in SO(3), \hat{p}(0), \hat{v}(0),\eta(0) \in \mathbb{R}^3$, $k_R, k_p, k_v$ are strictly positive scalar gains, and the innovation terms $\sigma_R$ and $y$ are given by
	\begin{align}
	&\sigma_R = \frac{1}{2}\sum_{i=1}^N k_i (p_i - p_c)^\times ( p_i - \hat{p} - \hat{R}y_i)  \label{eqn:sigmaR} \\
	& y    =   \sum_{i=1}^N k_i( p_i - \hat{p} - \hat{R}y_i) \label{eqn:sigmap}
	\end{align}
	with $k_i>0, \forall  i=1,2,\cdots,N$,  the landmark measurements $y_i$ in (\ref{eqn:output_y}), and  $p_c$   denoting the weighted center of the landmarks.
	Contrary to the work in \cite{BerkaneCDC2017,barrau2017invariant,berkane2019attitude}, where the attitude is updated intermittently, our attitude is updated continuously thanks to the auxiliary variable $\eta$ which takes care of the jumps upon the arrival of the landmark measurements.
	The introduction of the weighted center of the landmarks $p_c$ in $y$ and the dynamics of $\hat{p}$ allow us to decouple the position error dynamics and the attitude error dynamics, which is motivated by \cite{wang2020hybrid}. 
	
	Define the geometric estimation errors:  $\tilde{R}=R\hat{R}\T, \tilde{v} = v - \tilde{R}\hat{v}$ and $\tilde{p} = p -\tilde{R} \hat{p}- (I_3- \tilde{R})p_c$. Then, from the definitions of the output $y_i$ and the matrix $M$, the innovation terms $y$   and $\sigma_R$  can be rewritten in terms of the estimation errors as
	\begin{align}
	\sigma_R
	&=  -\frac{1}{2}\sum_{i=1}^N k_i (p_i - p_c)^\times  \tilde{R}\T (p_i-p_c)   =
	\psi(M\tilde{R})  \label{eqn:property-1} 	\\
	y &=    \sum_{i=1}^N k_i \tilde{R}\T \left( p -\tilde{R} \hat{p}- (I_3- \tilde{R})p_i\right) =    \tilde{R}\T \tilde{p} \label{eqn:property-2}
	\end{align}
	where we made use of  the facts: $\sum_{i=1}^N k_i (p_i - p_c)=0_{3\times 1}, x^\times x=0_{3\times 1}, \forall x\in \mathbb{R}^3$ and $x \times y = -y \times x= 2\psi(yx\T), \forall x,y\in \mathbb{R}^3 $. Note that when the weighted center of landmarks $p_c$ is located at the origin (\ie, $p_c = 0_{3\times 1}$), one has the traditional geometric position estimation error $\tilde{p}=p-\tilde{R}\hat{p}$. The advantage of our position estimation error is that the innovation term $y$ can be directly written in terms of the position estimation error. Moreover,  the introduction of $p_c$ in the expression of $\sigma_R$ given in (\ref{eqn:sigmaR})  results in $\sigma_R$ being only dependent on the attitude estimation error $\tilde{R}$  and not on   the position estimation error $\tilde{p}$ as shown in (\ref{eqn:property-1}).	
	
	In view of (\ref{eqn:R})-(\ref{eqn:v}),  (\ref{eqn:observer1}) and (\ref{eqn:property-1})-(\ref{eqn:property-2}), one has the following hybrid closed-loop system
	\begin{align}
	\underbrace{
		\begin{array}{ll}
		\dot{\tilde{R}} &= \tilde{R}(- \eta )^\times \\
		\dot{\eta}  &= 0_{3\times 1}    \\
		\dot{\tilde{p}} &=   \tilde{v}   \\
		\dot{\tilde{v}} &=   (I-\tilde{R})g
		\end{array}~}_{\tau\in[0, T_M]}
	\quad
	\underbrace{
		\begin{array}{ll}
		\tilde{R}^+ &= \tilde{R}   \\
		\eta^+    & = k_R\psi(M\tilde{R})  \\
		\tilde{p}^+ &=  \tilde{p}  -   k_p \tilde{p}    \\
		\tilde{v}^+ &=   \tilde{v} -   k_v \tilde{p}
		\end{array}~  }_{\tau \in \{0\} }
	\label{eqn:error1}
	\end{align}	
	Consider the new variable  $\mathsf{x} = [\tilde{p}\T,\tilde{v}\T]\T \in \mathbb{R}^6$, whose dynamics are given by
	\begin{equation}
	\begin{cases}
	\dot{\mathsf{x}}  =A \mathsf{x} + \delta_g & \tau\in[0, T_M] \\
	\mathsf{x}^+ = (I-KC) \mathsf{x} & \tau \in \{0\}
	\end{cases}  \label{eqn:hx1}
	\end{equation}
	where $\delta_g  := [0_{1\times 3},~ g\T (I_3-\tilde{R}) ]\T $, $K :=[ k_p I_3,
	k_v I_3]\T$, and the matrices $A$ and $C$ are given by
	\begin{align}
	A &= \begin{bmatrix}
	0_{3\times 3} & I_3  \\
	0_{3\times 3} & 0_{3\times 3}
	\end{bmatrix} ,
	C =    \begin{bmatrix}
	I_3~~ 0_{3\times 3}
	\end{bmatrix} . \label{eqn:AKC1}
	\end{align}
	The dynamics of $\mathsf{x}$ in (\ref{eqn:hx1}) can be seen as a linear hybrid system with a perturbation term $\delta_g$ induced by the gravity. Note that from the definition of $\delta_g$, this additional term will vanish as the attitude estimation error converges to $I_3$. Moreover, one can easily verify that the pair $(A,C)$ given in (\ref{eqn:AKC1}) is uniformly observable.

	Define the extended space $\mho: = SO(3)\times \mathbb{R}^3$. Then, we introduce the new state $x_1 = (\tilde{R},\eta,\mathsf{x},\tau)\in \mho \times \mathbb{R}^6 \times [0,T_M]$. From (\ref{eqn:error1}) and (\ref{eqn:hx1}), one obtains the hybrid closed-loop system
	$\mathcal{H}_1=(F_1,G_1,\mathcal{F}_1,\mathcal{J}_1)$ as follows:
	\begin{equation}
	\mathcal{H}_1:\begin{cases}
	\dot{x}_1 ~= F_1( {x}_1) &  x_1\in \mathcal{F}_1 \\
	x_1^+ \in G_1( {x}_1) & x_1 \in \mathcal{J}_1
	\end{cases}
	\label{eqn:closed-loop1}
	\end{equation}
	with $\mathcal{F}_1:= \mho \times \mathbb{R}^6 \times  [0, T_M] $, $\mathcal{J}_1:=\mho \times \mathbb{R}^6\times \{ 0\} $, and the flow and jump maps defined as
	\begin{align}
	F_1(x_1)&=\left(
	\tilde{R}(-\eta)^\times,
	0_{3\times 1},
	A\mathsf{x} + \delta_g ,
	-1
	\right) \label{eqn:F1}\\
	G_1(x_1)&= \left(
	\tilde{R},
	k_R\psi(M\tilde{R}),
	(I-KC) \mathsf{x},
	[T_m, T_M]
	\right)  \label{eqn:G1} .
	\end{align}
	Note that the flow set $\mathcal{F}_1$ and jump set $\mathcal{J}_1$ of $\mathcal{H}_1$ are closed, and $\mathcal{F}_1\cup \mathcal{J}_1 = \mho\times \mathbb{R}^6 \times [0,T_M] $. Moreover, with the introduction of the virtual timer $\tau$, the hybrid system $\mathcal{H}_1$ is autonomous and satisfies the hybrid basic conditions of \cite{goebel2009hybrid}.

	Define $\bar{M}: = \frac{1}{2}(\tr(M)I_3-M)$  with $M = \sum_{i=1}^N k_i (p_i-p_c)(p_i-p_c)\T$, which can be easily verified to be positive definite. Let us introduce a constant scalar associated to the minimum and maximum eigenvalues of the matrix $\bar{M}$ as $\varsigma_{M} := \lambda_{m}^{\bar{M}}/\lambda_{M}^{\bar{M}}>0$. Define the   closed set
	$
	\mathcal{A}: = \{x_1=(\tilde{R},\eta,\mathsf{x},\tau)\in \mho \times \mathbb{R}^6 \times [0,T_M]~|~
	\tilde{R}=I_3, \|\eta\|=0, \|\mathsf{x}\|=0   \}
	$.
	Now, one can state the following result:
	\begin{thm} \label{theo:1}
		Consider the hybrid dynamical system \eqref{eqn:closed-loop1}-\eqref{eqn:G1}.  Suppose that Assumption \ref{assum:1} - \ref{assum:2} hold, and there exists a symmetric positive definite  matrix $P$ satisfying
		\begin{equation}
		\Xi_P(\tau) :=A_g\T\Phi(\tau)\T P \Phi(\tau) A_g - P <0  \label{eqn:defXi_P}
		\end{equation}		
		for all $\tau \in [T_m,T_M] $ with $\Phi(\tau) = \exp(A \tau), A_g = (I-KC)$, $K=[k_p I_3, k_v I_3]\T$, and matrices $A$ and $C$ given in \eqref{eqn:AKC1}. Then, for any $0<\epsilon < 1 $, there exist constants $\eta^*,k_R^*>0$, such that for any  $|\tilde{R}(0)|_I  \leq \epsilon \sqrt{\varsigma_{M}} $, $\|\eta(0)\|\leq \eta^*$, $\mathsf{x}(0)\in \mathbb{R}^6$ and $k_R<  k_R^* $ the set $ {\mathcal{A}}$ is exponentially stable.
	\end{thm}

	\begin{pf}
		See   Appendix \ref{sec:theo1}.
	\end{pf}
	
	\begin{rem}
		To increase the basin of attraction for the attitude estimation error, one can choose $\|\eta(0)\|=0$ and $\varsigma_{M}=1$ (\ie, $M= k I_3$ with some constant $k>0$) through a proper construction of the matrix $M$, see for instance \cite{tayebi2013inertial}.
	\end{rem}
	\begin{rem}
		Note that a necessary condition for the existence of a symmetric positive definite  matrix $P$ satisfying \eqref{eqn:defXi_P} for all $ \tau \in [T_m,T_M]$ is that the pair $(\Phi(\tau),C)$ is observable for every $\tau \in [T_m,T_M]$, see  \cite{ferrante2016state}. The observability of the pair $(\Phi(\tau),C)$ for every $\tau \in [T_m,T_M]$   can be easily verified using $\Phi(\tau)=\exp(A \tau)$ and $A,C$ given in (\ref{eqn:AKC1}).  However, it is not straightforward to determine a sufficient condition for the existence of a solution for the optimization problem $\Xi_P(\tau)< 0, \forall \tau \in [T_m,T_M] $.
		This optimization problem can be solved using the polytopic embedding technique proposed in \cite{ferrante2016state} and the finite-dimensional LMI approach proposed in the recent work \cite{sferlazza2019time}. A complete procedure for solving this infinite-dimensional optimization problem, adapted from the  work in \cite{sferlazza2019time}, is provided in Appendix \ref{sec:solvingXoP}. In practice, one can start with the design of a  gain $K$ such that the eigenvalues of $\Phi(\tau) (I-KC)$ are inside the unit circle for some $\tau \in [T_m,T_M]$, and apply the procedure in Appendix \ref{sec:solvingXoP} with this gain $K$ to check whether \eqref{eqn:defXi_P} is satisfied for all $\tau \in [T_m,T_M]$.  
	\end{rem}
	
	\subsubsection*{Unknown gravity case (with gravity estimation):}
	In some applications, the gravity vector $g$ may not be available or not accurately known. Hence, it is of great interest to design observers without the knowledge of the gravity vector. To solve this problem, a new hybrid observer with gravity vector estimation is proposed.	Let $\hat{g} \in \mathbb{R}^3$ be the estimate of the gravity vector $g$, and $\tilde{g}: = g - \tilde{R}\hat{g}$ be the gravity vector estimation error. We propose the following hybrid observer:
	\begin{align}
	\underbrace{
		\begin{array}{ll}
		\dot{\hat{R}} &= \hat{R}(\omega + \hat{R}\T \eta )^\times \\
		\dot{\eta}  &= 0_{3\times 1}    \\
		\dot{\hat{p}} &= \eta^\times (\hat{p}-p_c) + \hat{v}   \\
		\dot{\hat{v}} &=  \eta^\times \hat{v} +\hat{g} + \hat{R}a   \\
		\dot{\hat{g }} & = \eta^\times \hat{g }
		\end{array}~}_{\tau\in[0, T_M]}
	\quad
	\underbrace{
		\begin{array}{ll}
		\hat{R}^+ &= \hat{R}   \\
		\eta^+    & = k_R \sigma_R  \\
		\hat{p}^+ &=  \hat{p}  +  k_p   y     \\
		\hat{v}^+ &=    \hat{v}  +  k_v  y \\
		\hat{g }^+ &=   \hat{g }  + k_g  y
		\end{array}~  }_{\tau \in \{0\} }
	\label{eqn:observer2}
	\end{align}
	where $\hat{R}(0)\in SO(3), \hat{p}(0), \hat{v}(0), \eta(0), \hat{g}(0) \in \mathbb{R}^3$, the constant scalar gains $k_R, k_p, k_v, k_g >0$, and the innovation terms $\sigma_R$ and $y$ are given in
	\eqref{eqn:sigmaR}  and \eqref{eqn:sigmap}, respectively. Compared to the hybrid observer (\ref{eqn:observer1}), in the dynamics  of $\hat{v}$ the gravity vector $g$ has been replaced by $\hat{g}$ which is obtained from an appropriately designed adaptation law. From (\ref{eqn:observer2}), the estimated gravity vector $\hat{g}$ is continuously updated using $\eta$ when the landmark measurements are not available (\textit{i.e.,} in the flow set), and discretely updated using $y$ upon the arrival of the landmark measurements (\textit{i.e.,} in the jump set).
	
	From (\ref{eqn:observer1}) and (\ref{eqn:observer2}), one obtains the same hybrid closed-loop dynamics for $\tilde{R}$ and $\eta$ as (\ref{eqn:error1}). Let us introduce a new variable $\mathsf{x}:=[\tilde{p}\T, \tilde{v}\T, \tilde{g }\T]\T \in \mathbb{R}^9$. From (\ref{eqn:R})-(\ref{eqn:v}),  (\ref{eqn:observer2}),  the hybrid dynamics of $\mathsf{x}$ are given as follows
	\begin{equation}
	\begin{cases}
	\dot{\mathsf{x}}  =A \mathsf{x}  & \tau\in[0, T_M] \\
	\mathsf{x}^+ = (I-KC) \mathsf{x} & \tau \in \{0\}
	\end{cases}  \label{eqn:hx2}
	\end{equation}
	where   $K :=[ k_p I_3,
	k_v I_3, k_g I_3]\T$, and the matrices $A$ and  $C$ are given by
	\begin{align}
	&A = \begin{bmatrix}
	0_{3\times 3} & I_3 & 0_{3\times 3}\\
	0_{3\times 3} & 0_{3\times 3} & I_3\\
	0_{3\times 3} & 0_{3\times 3}  & 0_{3\times 3}
	\end{bmatrix} , C=   \begin{bmatrix}
	I_3\\
	0_{3\times 3}\\	
	0_{3\times 3}
	\end{bmatrix}\T .
	\label{eqn:AKC2}
	\end{align}
	The difference between (\ref{eqn:hx1}) and (\ref{eqn:hx2}) is that $\delta_{g}=0_{9\times 1}$ and $\mathsf{x}\in \mathbb{R}^9$ in (\ref{eqn:hx2}).
	Hence, in view of  (\ref{eqn:hx1})  and (\ref{eqn:hx2}), one obtains the same form of hybrid closed-loop system as  $\mathcal{H}_1$ in (\ref{eqn:closed-loop1}) for observer \eqref{eqn:observer2}  except that $\delta_{g}=0_{9\times 1}$ and $\mathsf{x}\in \mathbb{R}^9$. Let us introduce the following closed set:
	$
	\bar{\mathcal{A}}: = \{x_1=(\tilde{R},\eta,\mathsf{x},\tau)\in \mho \times \mathbb{R}^9 \times [0,T_M]~| ~
	\tilde{R}=I_3, \|\eta\|= \|\mathsf{x}\|= \|\tilde{g}\|=0  \}
	$.	Now, one can state the following result:
	\begin{thm} \label{theo:6}
		Consider the hybrid observer \eqref{eqn:observer2} for the system \eqref{eqn:R}-\eqref{eqn:v}.  Suppose that  Assumption \ref{assum:1} -\ref{assum:2} hold, and there exists a symmetric positive definite  matrix $P$ satisfying \eqref{eqn:defXi_P} for all $\tau \in [T_m,T_M] $ with $\Phi(\tau) = \exp(A \tau)$, $K=[k_pI_3, k_vI_3, k_gI_3]\T$, and  matrices $A$ and $C$ given in  \eqref{eqn:AKC2}. Then, for any $0<\epsilon < 1 $, there exist constants $\eta^*,k_R^*>0$, such that for any  $|\tilde{R}(0)|_I  \leq \epsilon \sqrt{\varsigma_{M}} $, $\|\eta(0)\|\leq \eta^*$, $\mathsf{x}(0)\in \mathbb{R}^9$ and $k_R<  k_R^* $ the set $\bar{\mathcal{A}}$ is exponentially stable.
	\end{thm} 
	The proof of Theorem \ref{theo:6} can be conducted using the same steps as in the proof of Theorem \ref{theo:1}, and is therefore omitted here.

	\subsection{Variable-gain Design} \label{sec:observerII}
	\subsubsection*{Known gravity case (without gravity estimation):}
	In the previous subsection, the fixed-gain design approach was considered based on an infinite-dimensional optimization. However, in practice, it is not straightforward to tune the gain parameters to satisfy condition (\ref{eqn:defXi_P}) for all $\tau \in [T_m,T_M]$. Moreover, since this condition is dependent on the parameters $T_m,T_M$, one may need to redesign the gains every time these parameters change.  In this subsection, we a propose different hybrid observer  relying on variable gains automatically designed via continuous-discrete Riccati equations.
	
	We propose the following hybrid nonlinear observer:
	\begin{align}
	\underbrace{
		\begin{array}{ll}
		\dot{\hat{R}} &= \hat{R}(\omega + \hat{R}\T \eta )^\times   \\
		\dot{\eta}  &= 0_{3\times 1}  \\
		\dot{\hat{p}} &= \eta^\times (\hat{p}-p_c) + \hat{v}   \\
		\dot{\hat{v}} &= \eta^\times  \hat{v} + g  + \hat{R}a
		\end{array}~}_{\tau \in [T_m, T_M] }
	\quad
	\underbrace{
		\begin{array}{ll}
		\hat{R}^+ &= \hat{R}  \\
		\eta^+    & = k_R\sigma_R   \\
		\hat{p}^+ &=  \hat{p}  + \hat{R} K_p \hat{R}\T y     \\
		\hat{v}^+ &=  \hat{v}  + \hat{R}  K_v \hat{R}\T  y
		\end{array}~  }_{\tau \in \{0\}}
	\label{eqn:observer3}
	\end{align}
	where $\hat{R}(0)\in SO(3), \hat{p}(0), \hat{v}(0),\eta(0)  \in \mathbb{R}^3$ and $k_R>0$. The innovation terms $\sigma_R$ and $y$ are given in (\ref{eqn:sigmaR}) and (\ref{eqn:sigmap}), respectively. The main difference with respect to observer (\ref{eqn:observer1}) is that the gain matrices $K_p, K_v \in \mathbb{R}^{3\times 3}$ in \eqref{eqn:observer3} are time-varying matrices.  In view of (\ref{eqn:R})-(\ref{eqn:v})  and  (\ref{eqn:observer3}), one has the following hybrid closed-loop system:
	\begin{align}
	\underbrace{
		\begin{array}{ll}
		\dot{\tilde{R}} &= \tilde{R}(- \eta )^\times \\
		\dot{\eta}  &= 0_{3\times 1}    \\
		\dot{\tilde{p}} &=   \tilde{v}   \\
		\dot{\tilde{v}} &=   (I-\tilde{R})g
		\end{array}~}_{\tau\in[0, T_M]}
	\quad
	\underbrace{
		\begin{array}{ll}
		\tilde{R}^+ &= \tilde{R}   \\
		\eta^+    & = k_R\psi(M\tilde{R})  \\
		\tilde{p}^+ &=  \tilde{p}  -   R K_p R\T \tilde{p}    \\
		\tilde{v}^+ &=   \tilde{v} -   R K_v R\T \tilde{p}
		\end{array}~  }_{\tau \in \{0\} }
	\label{eqn:error2}
	\end{align}

	Let us introduce a new variable $ {\mathsf{x}}  = [\tilde{p}\T R,\tilde{v}\T R]\T \in \mathbb{R}^6$.
	From (\ref{eqn:R}) and (\ref{eqn:error2}), the hybrid dynamics of $\mathsf{x} $ are given by
	\begin{equation}
	\begin{cases}
	\dot{\mathsf{x}}   =A_t \mathsf{x} + \bar{\delta}_{g}  & \tau \in [0,T_M]\\
	\mathsf{x}^+  = (I-KC) \mathsf{x} & \tau \in \{0\}
	\end{cases}  \label{eqn:hx3}
	\end{equation}
	where $\bar{\delta}_g  := [0_{1\times 3}~ g\T (R-\hat{R}) ]\T$, $K:=[K_p\T,K_v\T]\T$, and the matrices $A_t$ and $C$ are given by
	\begin{align}
	A_t &= \begin{bmatrix}
	-\omega^\times & I_3  \\
	0_{3\times 3} & -\omega^\times
	\end{bmatrix} ,
	C  = \begin{bmatrix}
	I_3 &
	0_{3\times 3}
	\end{bmatrix} . \label{eqn:AKC3}
	\end{align}
	Similar to (\ref{eqn:hx1}), the dynamics of $\mathsf{x}$ in (\ref{eqn:hx3}) can be seen as a linear hybrid system with an additional perturbation term $\bar{\delta}_{g}$ induced by the gravity. One can also show that $\bar{\delta}_{g}$ vanishes as the attitude estimation error converges to $I_3$. The main difference between (\ref{eqn:hx1}) and (\ref{eqn:hx3}) is that matrix $A_t$ in (\ref{eqn:hx3}) is time-varying.  Let us design the gain matrix $K$ as
	\begin{equation}
	K = PC\T(CPC\T + Q_t)^{-1} \label{eqn:Kt2}
	\end{equation} 	
	where $P$ is the solution of the following continuous-discrete Riccati equation
	\begin{subequations}
		\begin{align}
		\dot{P} ~~&= A_tP + PA_t\T + V_t, \quad  \tau\in [0,T_M] \label{eqn:CDRE-C}\\
		P^+ &= P-PC\T(CPC\T +Q_t)^{-1}C P,~ \tau \in \{0\} \label{eqn:CDRE-D}
		\end{align}
	\end{subequations}
	where $P(0)$ is positive definite, $V_t$ is continuous, $Q_t, V_t$ are uniformly positive definite, and the matrices $(A_t, C)$ are given by (\ref{eqn:AKC3}). The following lemma, adapted from \cite{deyst1968conditions,barrau2017invariant}, provides sufficient conditions for the existence of the solution of the continuous-discrete Riccati equations (\ref{eqn:CDRE-C})-(\ref{eqn:CDRE-D}).
	\begin{lem}{\cite{deyst1968conditions,barrau2017invariant}} \label{lem:lemma-CDRE}
		Consider the pair $(A_t,C)$, and let matrix $A_t, V_t$ be continuous and matrices $V_t$ and $Q_t$ be uniformly positive definite. If there exist constants $\mu_{\phi}, {\mu}_v,  {\mu}_V, {\mu}_q,  {\mu}_Q \in \mathbb{R}_{>0}$ and $ \Gamma \in \mathbb{N}_{>0}$ such that for  all $j \geq \Gamma$
		\begin{subequations}
			\begin{align}
			&\mu_{\phi} I \leq (\Phi^{t_{j+1}}_{t_j})\T \Phi^{t_{j+1}}_{t_j} \label{eqn:CDRE_condition1}\\
			&{\mu}_v I \leq  \int_{t_{j -  \Gamma} }^{t_{j}} \Phi^{t_j}_s {V}_t(s)(\Phi^{t_j}_s)\T ds \leq  {\mu}_V I \label{eqn:CDRE_condition2}\\ 
			&\mu_q I \leq  \sum_{i=j-\Gamma}^{j}  (\Phi^{t_i}_{t_{j}})\T C \T Q_t^{-1} (t_i) C  \Phi^{t_i}_{t_{j}}  \leq \mu_Q I \label{eqn:CDRE_condition3}
			\end{align}
		\end{subequations}
		where  $\Phi^t_{s} = \Phi(t,s)$ denotes the square matrix defined by $\Phi^t_{t} = I, \partial (\Phi^t_s)/\partial t = A_t\Phi^t_s$. Then, the solution $P$ to  (\ref{eqn:CDRE-C})-(\ref{eqn:CDRE-D}) exists, and there exist constants $0< p_m \leq p_M < \infty $  such that $p_m I \leq P \leq p_M I$ for all $ t \geq 0$. 
	\end{lem} 
	\begin{rem}	
		As pointed out in \cite{jazwinski1970}, a continuous-discrete filter can be embedded in a discrete filter. Indeed, equations  (\ref{eqn:CDRE-C})-(\ref{eqn:CDRE-D}) can be cast into a discrete setting as shown in \cite[Section 7.2]{jazwinski1970}, which allows to use the results of \cite{deyst1968conditions} with slight modifications to retrieve the conditions in Lemma \ref{lem:lemma-CDRE} and \cite[Theorem 3]{barrau2017invariant}. 
	\end{rem}
	Define the new state $x_2 = (\tilde{R},\eta, {\mathsf{x}},\tau)$. From (\ref{eqn:error2}) and (\ref{eqn:hx3}), one obtains the following hybrid closed-loop system:
	$\mathcal{H}_2=(F_2,G_2,\mathcal{F}_2,\mathcal{J}_2)$:
	\begin{equation}
	\mathcal{H}_2:\begin{cases}
	\dot{x}_2 ~= F_2( {x}_2) &  x_2\in \mathcal{F}_2\\
	x_2^+ \in G_2({x}_2) & x_2\in \mathcal{J}_2
	\end{cases}
	\label{eqn:closed-loop2}
	\end{equation}
	with $\mathcal{F}_2:= \mho \times \mathbb{R}^6 \times  [0, T_M] $, $
	\mathcal{J}_2:= \mho\times \mathbb{R}^6\times \{ 0\} $, and the following flow and jump maps:
	\begin{align}
	F_2(x_2)&=(
	\tilde{R}(-\eta)^\times,
	0_{3\times 1},
	A_t\mathsf{x} + \bar{\delta}_g ,
	-1
	)  \label{eqn:F2}\\
	G_2(x_2)&= (
	\tilde{R},
	k_R\psi(M\tilde{R}),
	(I-KC) \mathsf{x},
	[T_m, T_M]
	). \label{eqn:G2}
	\end{align}
	Similar to (\ref{eqn:closed-loop1}), one can show that the hybrid system $\mathcal{H}_2$ satisfies the hybrid basic conditions.	
	Now, one can  state   the following result:
	\begin{thm} \label{theo:2}
		Consider the hybrid dynamical system \eqref{eqn:closed-loop2}-\eqref{eqn:G2}.  Let Assumption  \ref{assum:1} - \ref{assum:2} hold. Assume that $\omega(t)$ and $V_t$ are continuous and bounded, $Q_t$ is bounded, and $V_t,~Q_t$  are uniformly positive definite.
		Then, for any $0<\epsilon < 1 $, there exist constants $\eta^*,k_R^*>0$, such that for any  $|\tilde{R}(0)|_I  \leq \epsilon \sqrt{\varsigma_{M}} $, $\|\eta(0)\|\leq \eta^*$, $\mathsf{x}(0)\in \mathbb{R}^6$ and $k_R<  k_R^* $ the set $ {\mathcal{A}}$ is exponentially stable.
	\end{thm}
	\begin{pf}
		See  Appendix \ref{sec:theo2}.
	\end{pf}
	\begin{rem}\label{rem:gaintuning}
		It is important to mention that the proposed observer \eqref{eqn:observer3} is deterministic 
		and the results in Theorem \ref{theo:2} hold for any $Q_t$ and $V_t$ uniformly positive definite. However, in practice, the matrices $V_t$ and $Q_t$ can be tuned using the (approximately known) covariance of the state and output noise as in the Kalman filter. Let $\omega$ and $a$ be the noisy measurements, and replace $\omega$ by $\omega+n_{\omega}$ in \eqref{eqn:R} and $a$ by $a+n_a$ in \eqref{eqn:v} with $n_\omega, n_a \in \mathbb{R}^3$ denoting the noise signals on the angular velocity and acceleration measurements. Differentiating $R\T \tilde{p}$ and $R\T \tilde{v}$, in view of (\ref{eqn:R})-(\ref{eqn:v})  and  (\ref{eqn:observer3}), and linearizing around $\mathsf{x}= 0$, one obtains 
		\begin{align*}
		\dot{\mathsf{x}}  ={A}_t \mathsf{x} + \bar{\delta}_g + G_t \begin{bmatrix}
		n_\omega\\
		n_a
		\end{bmatrix}, G_t   =  \begin{bmatrix}
		(\hat{R}\T(\hat{p}-p_c)) ^\times & 0_{3\times 3} \\
		(\hat{R}\T \hat{v})^\times & I_3
		\end{bmatrix}
		\end{align*}
		with $p_c = \sum_{i=1}^{N} k_i p_i$ and $A_t$ defined in \eqref{eqn:AKC3}.
		Replacing $y_i=R\T(p_i-p)$ by $y_i = R\T(p_i-p) + n_y^i$ in \eqref{eqn:output_y} for all $i=1,2,\dots,N$ with $n_y^i \in \mathbb{R}^3$ denoting the noise in the landmark measurements, the virtual output $y$ in \eqref{eqn:sigmap} can be rewritten as
		$
		y  = C \mathsf{x} - \sum_{i=1}^N k_i \hat{R} n_y^i,   
		$
		with $C$ defined in \eqref{eqn:AKC3}.
		Hence, the matrices $V_t$ and $Q_t$ can be chosen as
		\begin{align}
		V_t = G_t \text{Cov}\left( \begin{bmatrix}
		n_\omega\\
		n_a
		\end{bmatrix} \right)  G_t\T, 
		Q_t  = \sum_{i=1}^N k_i^2 \hat{R} \text{Cov}\left( n_y^i  \right)\hat{R}\T \label{eqn:VQ1}
		\end{align}
		where we made the assumption that the noise signals in the landmark measurements are uncorrelated. In practice, a small positive definite matrix can be added to $V_t$ in \eqref{eqn:VQ1} to ensure that $V_t$ is uniformly positive definite.
	\end{rem}
	
	\subsubsection*{Unknown gravity case (with gravity estimation):}
	On the other hand, to handle the unknown gravity case, we propose the following new hybrid nonlinear observer with gravity vector estimation:
	\begin{align}
	\underbrace{
		\begin{array}{ll}
		\dot{\hat{R}} &= \hat{R}(\omega + \hat{R}\T \eta )^\times \\
		\dot{\eta}  &= 0_{3\times 1}    \\
		\dot{\hat{p}} &= \eta^\times (\hat{p}-p_c) + \hat{v}   \\
		\dot{\hat{v}} &=  \eta^\times \hat{v} +\hat{g} + \hat{R}a   \\
		\dot{\hat{g }} & = \eta^\times \hat{g }
		\end{array}~}_{\tau\in[0, T_M]}
	\quad
	\underbrace{
		\begin{array}{ll}
		\hat{R}^+ &= \hat{R}   \\
		\eta^+    & = k_R \sigma_R  \\
		\hat{p}^+ &=  \hat{p}  + \hat{R}K_p \hat{R}\T y     \\
		\hat{v}^+ &=    \hat{v}  + \hat{R}K_v \hat{R}\T y \\
		\hat{g }^+ &=   \hat{g }  + \hat{R}K_g \hat{R}\T y
		\end{array}~  }_{\tau \in \{0\} }
	\label{eqn:observer4}
	\end{align}
	where $\hat{R}(0)\in SO(3), \hat{p}(0), \hat{v}(0), \hat{g}(0), \eta(0) \in \mathbb{R}^3$, $k_R>0$, and $K_p,K_v,K_{g}\in \mathbb{R}^{3\times 3}$ to be designed. The innovation terms $\sigma_R,y$ are given in (\ref{eqn:sigmaR}) and (\ref{eqn:sigmap}), respectively.
	
	Consider the new variable $ \mathsf{x} =[\tilde{p}\T R, \tilde{v}\T R,  \tilde{g }\T R]\T \in \mathbb{R}^9$.
	Then, from (\ref{eqn:R})-(\ref{eqn:v}) and (\ref{eqn:observer4}),  the dynamics of $\mathsf{x}$ are given by
	\begin{equation}
	\begin{cases}
	\dot{\mathsf{x}}   =A_t \mathsf{x}  & \tau \in [0,T_M]\\
	\mathsf{x}^+  = (I-KC) \mathsf{x} & \tau \in \{0\}
	\end{cases}  \label{eqn:hx4}
	\end{equation}
	where $K:=[K_p\T,K_v\T,K_g\T]\T$,  and the matrices $A_t$ and $C$ are given by
	\begin{align}
	A_t  = \begin{bmatrix}
	-\omega^\times & I_3 & 0_{3\times 3}\\
	0_{3\times 3} & -\omega^\times & I_3\\
	0_{3\times 3} & 0_{3\times 3}  & -\omega^\times
	\end{bmatrix},
	C = \begin{bmatrix}
	I_3 \\
	0_{3\times 3} \\
	0_{3\times 3}
	\end{bmatrix}\T.
	\label{eqn:AKC4}
	\end{align}
	The gain matrix $K$ is designed as
	\begin{equation}
	K  = P C\T(CPC\T +Q_t)^{-1} \label{eqn:Kt}
	\end{equation}
	where $P$ is the solution of the continuous-discrete Riccati equation (\ref{eqn:CDRE-C})-(\ref{eqn:CDRE-D}), with $P(0)$ positive definite,  $V_t, Q_t$ uniformly positive definite, and $(A_t, C)$ given by (\ref{eqn:AKC4}).
	In view of (\ref{eqn:hx3})  and (\ref{eqn:hx4}), one obtains a hybrid closed-loop system in the same form as $\mathcal{H}_2$ in (\ref{eqn:closed-loop2}) for observer \eqref{eqn:observer4}  except that $\bar{\delta}_{g}=0_{9\times 1}$ and $\mathsf{x}\in \mathbb{R}^9$. 
	Now, one can state the following result:
	\begin{thm} \label{theo:10}
		Consider the hybrid observer \eqref{eqn:observer4} for the system \eqref{eqn:R}-\eqref{eqn:v}.  Let  Assumption \ref{assum:1} - \ref{assum:2} hold, Assume that $\omega(t)$ and $V_t$ are continuous and bounded, $Q_t$ is bounded, and $V_t,~Q_t$  are uniformly positive definite.
		Then, for any $0<\epsilon < 1 $, there exist constants $\eta^*,k_R^*>0$, such that for any  $|\tilde{R}(0)|_I  \leq \epsilon \sqrt{\varsigma_{M}} $, $\|\eta(0)\|\leq \eta^*$, $\mathsf{x}(0)\in \mathbb{R}^9$ and $k_R<  k_R^* $ the set $\bar{\mathcal{A}}$ is exponentially stable.
	\end{thm}
	
	\begin{rem} 
		The proof of Theorem \ref{theo:10} can be conducted using the same steps as in the proof of Theorem \ref{theo:2}, and is therefore omitted here. Using the same steps as in Remark \ref{rem:gaintuning}, one can show that, in this case, $V_t$ and $Q_t$ can be related to the covariance matrices of the measurements noise as follows:
		\begin{align}
		V_t= G_t \text{Cov}\left( \begin{bmatrix}
		n_\omega\\
		n_a
		\end{bmatrix} \right)   G_t\T,  
		Q_t = \sum_{i=1}^N k_i^2 \hat{R} \text{Cov}\left( n_y^i  \right)\hat{R}\T   \label{eqn:VQ2}
		\end{align}
		with $n_\omega,n_a,n_y^i \in \mathbb{R}^3, i=1,2,\dots,N$ denoting the noise signals on the angular velocity, linear acceleration and landmark measurements, respectively, and
		\[
		G_t  =  \begin{bmatrix}
		(\hat{R}\T(\hat{p}-p_c))^\times & 0_{3\times 3} \\
		(\hat{R}\T \hat{v})^\times & I_3    \\
		(\hat{R}\T \hat{g})^\times & 0_{3\times 3}
		\end{bmatrix}. 
		\]
	\end{rem}

	\section{Simulation Results}\label{sec:simulation}
	
	In this section, simulation results are presented to illustrate the performance of the proposed hybrid observers. We refer to the  hybrid observer (\ref{eqn:observer1})   as `HINO1-F', the  hybrid observer (\ref{eqn:observer2})  as `HINO2-F',  the  hybrid observer (\ref{eqn:observer3})   as `HINO1-V', and the  hybrid observer (\ref{eqn:observer4})   as `HINO2-V'. Moreover, we refer to the invariant observer proposed in \cite{barrau2017invariant} as `IEKF'.
	
	We consider an autonomous vehicle moving on the `8'-shape trajectory given by $p(t)=10[\sin(t), \sin(t)\cos(t), 1]\T$ (m), with the initial rotation  $R(0)=I_3$ and the angular velocity $\omega(t) = [\sin(0.3\pi),0.1,\cos(0.3\pi)]\T$ (rad/s). The same initial conditions are considered for each observer $\hat{R}(0) = \mathcal{R}_a(0.1\pi,u)$ with $u\in \mathbb{S}^2$, $ \eta(0)=\hat{v}(0)=\hat{p}(0)=\hat{g}(0) = 0_{3\times 1}$, and $P(0)=I$. There are $N=25$ landmarks which are randomly selected on the ground such that Assumption \ref{assum:1} holds.  We consider continuous IMU measurements and intermittent landmark position measurements with $T_m = 0.04s$ and $T_M=0.06 s$ (about $20$Hz sampling rate). Moreover, additive white Gaussian noise has been considered  with $\text{Cov}(n_\omega) = 0.0001 I_3, \text{Cov}(n_a) = 0.01 I_3$ and $ \text{Cov}(n_y^i) = 0.01 I_3, \forall  i=1,2,\cdots N$ for the gryo, accelerometer and landmark measurements, respectively.  The same parameters  $k_i=1/N, i=1,\hdots,N$ and $k_R=1.2$ are chosen for each of the proposed observers. Moreover, for the fixed-gain observers,  we pick $k_p = 0.5, k_v = 1.0$ and $k_{g} = 0.6$, such that, for both observers, there exists a matrix $P$ satisfying $\Xi_P(\tau)<0, \forall \tau\in [T_m,T_M]$.
	For the variable-gain observers, matrices
	$V_t $ and $ Q_t $ are chosen as (\ref{eqn:VQ1}) for HINO1-V and as (\ref{eqn:VQ2})  for HINO2-V.  For the IEKF,  the gain matrices  are chosen using the same covariance of the measurements noise as per Section V.B in \cite{barrau2017invariant}.

\begin{figure}[!thpb]
		\centering
		\includegraphics[width=0.7\linewidth]{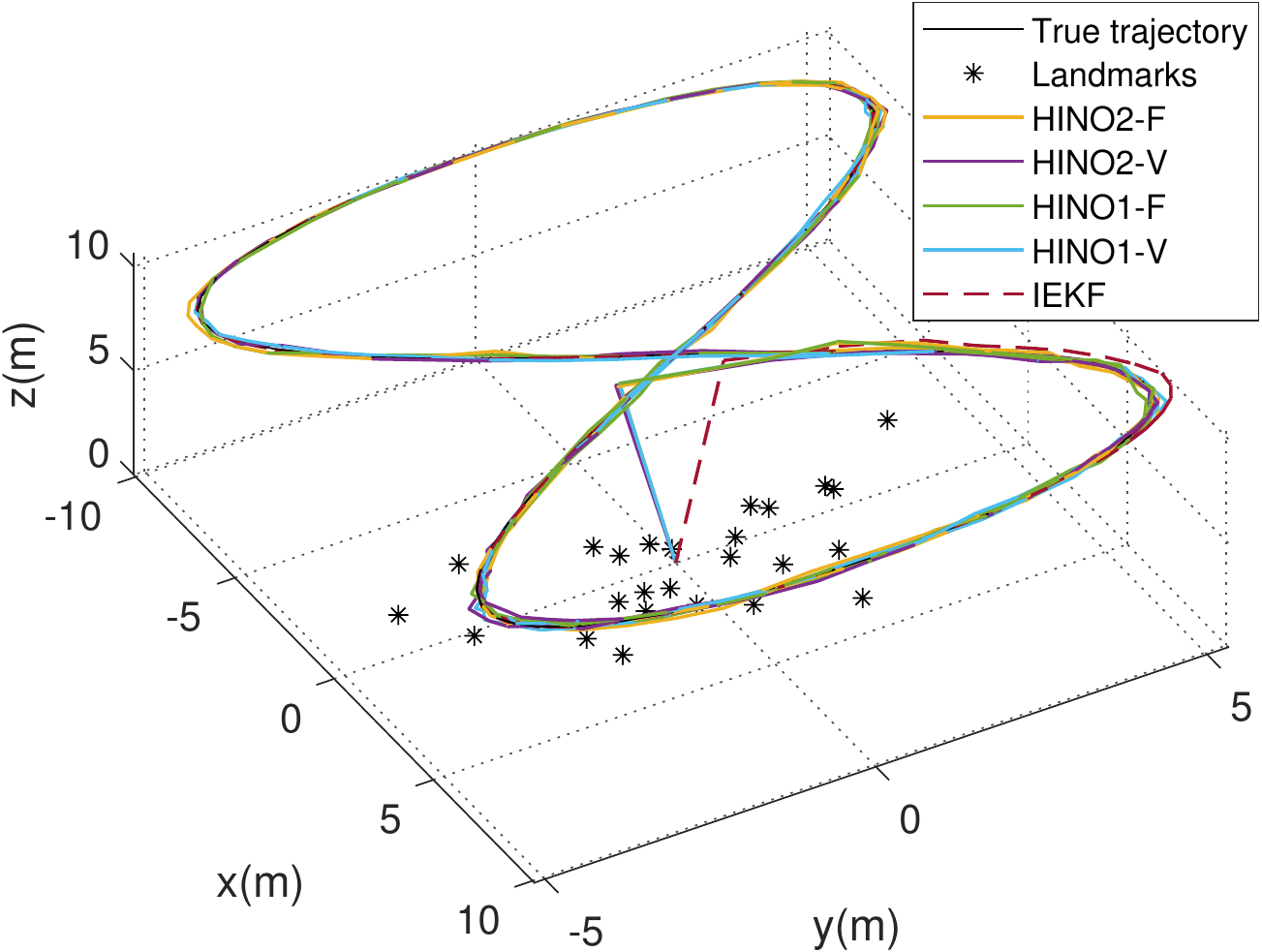} \\
		\includegraphics[width=0.95\linewidth]{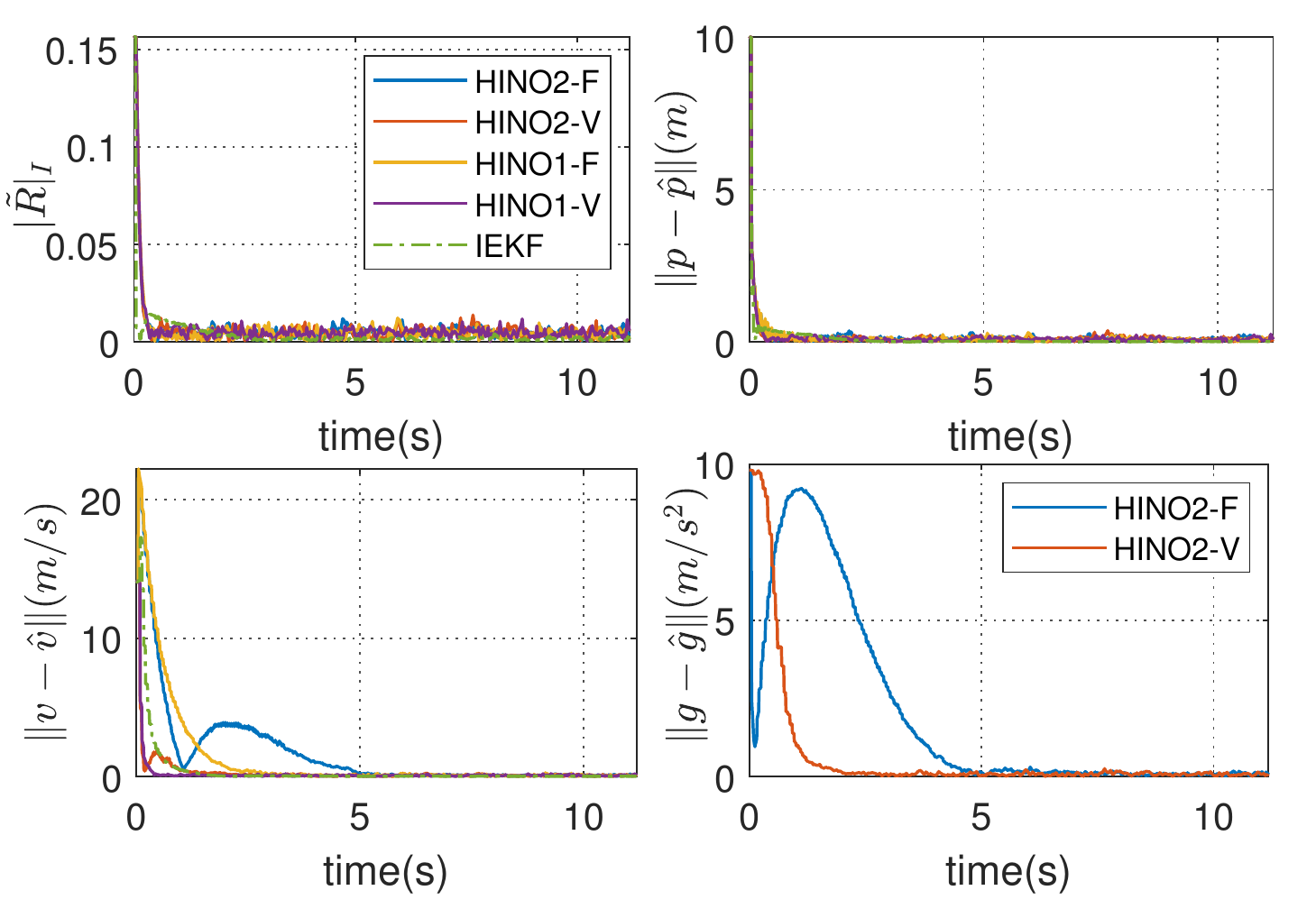}
		\caption{Simulation results of three-dimensional trajectories and estimation errors of rotation, position, linear velocity and gravity vector using intermittent landmark position measurements. }
		\label{fig:simulation}
	\end{figure}

	Simulation results are shown in Fig. \ref{fig:simulation}. As one can see, the estimated states from the proposed hybrid observers and IEKF converge, after a few seconds, to the vicinity of the real state. The execution time of each observer, using an Intel Core i7-3540M running at 3.00GHz, is given in the following table : \par
	\renewcommand{\arraystretch}{1.2}
	\begin{table}[!htpb]	\centering  \footnotesize	
		\scalebox{0.89}{\begin{tabular}{c|c|c|c|c|c}
				\hline
				$N$ & HINO1-F&HINO1-V  &HINO2-F  &HINO2-V  & IEKF  \\
				\hline
				25 & 0.0053s  &0.0084s  &0.0066s  &0.0095s  &0.0102s  \\
				\hline
				100 & 0.0061s  &0.0089s  &0.0069s  &0.0120s  &0.0138s  \\
				\hline
		\end{tabular}}
	\end{table}
	From the table, one can see that the computational costs of our fixed-gain observers are lower than the computational costs of our variable-gain observers and IEKF.
	This is due to the online computations required for solving the continuous-discrete Riccati equations. Moreover, the IEKF comes with the highest computational cost, which is mainly due to the computation of the inverse of a potentially high-dimensional matrix $S \in \mathbb{R}^{3N \times 3N}$ (see Eqn. (35) in \cite{barrau2017invariant}), especially when the number of landmarks $N$ is large.
	
	\section{Experimental Results}\label{sec:experimental}	
\begin{figure}[!thpb]
		\centering
		\includegraphics[width=0.85\linewidth]{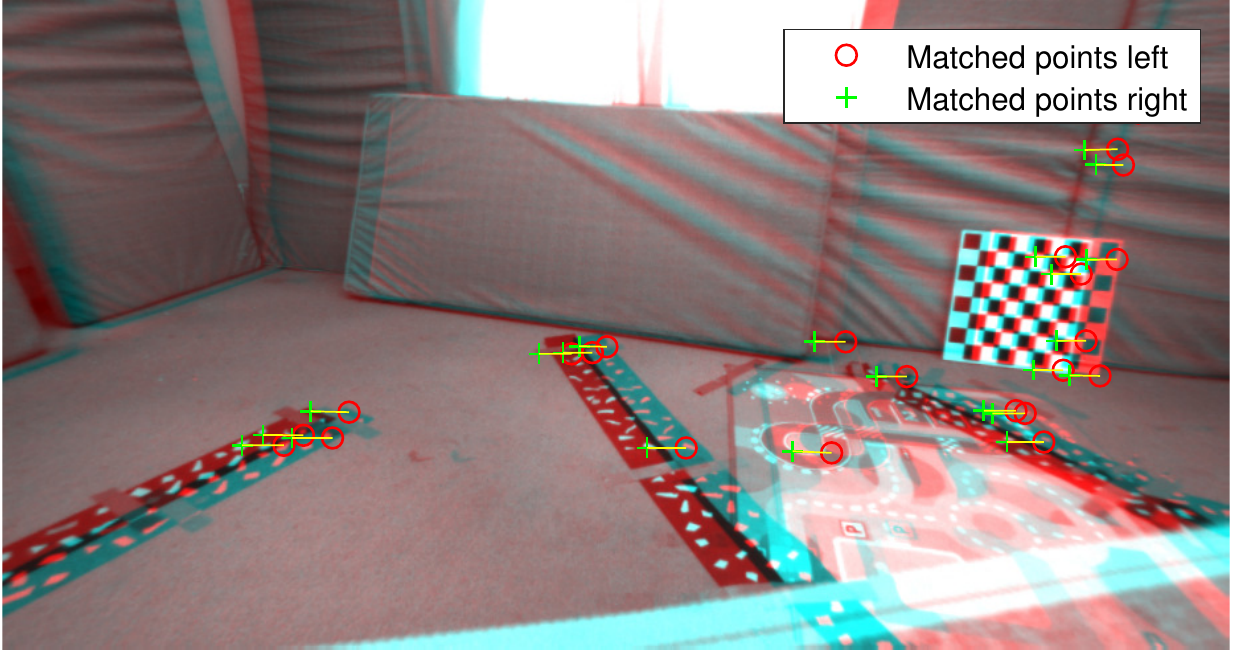}
		\caption{Example of features detection and tracking in the left and right images of a stereo camera using the Computer Vision System Toolbox. Pictures come from the EuRoc dataset\cite{Burri25012016}.}
		\label{fig:features}
	\end{figure}

	To further validate the performance of our proposed hybrid observers, we applied our algorithms to real data from the EuRoc dataset \cite{Burri25012016}, where the trajectories are generated by a real flight of a quadrotor. This dataset includes stereo images, IMU measurements and ground truth (obtained using Vicon motion capture system). The sampling rate of the IMU measurements from ADIS16448 is 200Hz and the sampling rate of the stereo images from MT9V034 is 20Hz. 
	The features are tracked via the Kanade-Lucas-Tomasi (KLT) tracker  \cite{shi1994good} using minimum eigenvalue feature detection (see Fig. \ref{fig:features}). More details about the EuRoc dataset can be found in \cite{Burri25012016} and details about the experimental setup for observers implementation can be found in \cite{wang2020hybrid}.

	\begin{figure*}[thpb]
		\centering 
		\begin{minipage}{0.92\linewidth}
				\includegraphics[width=0.40\linewidth]{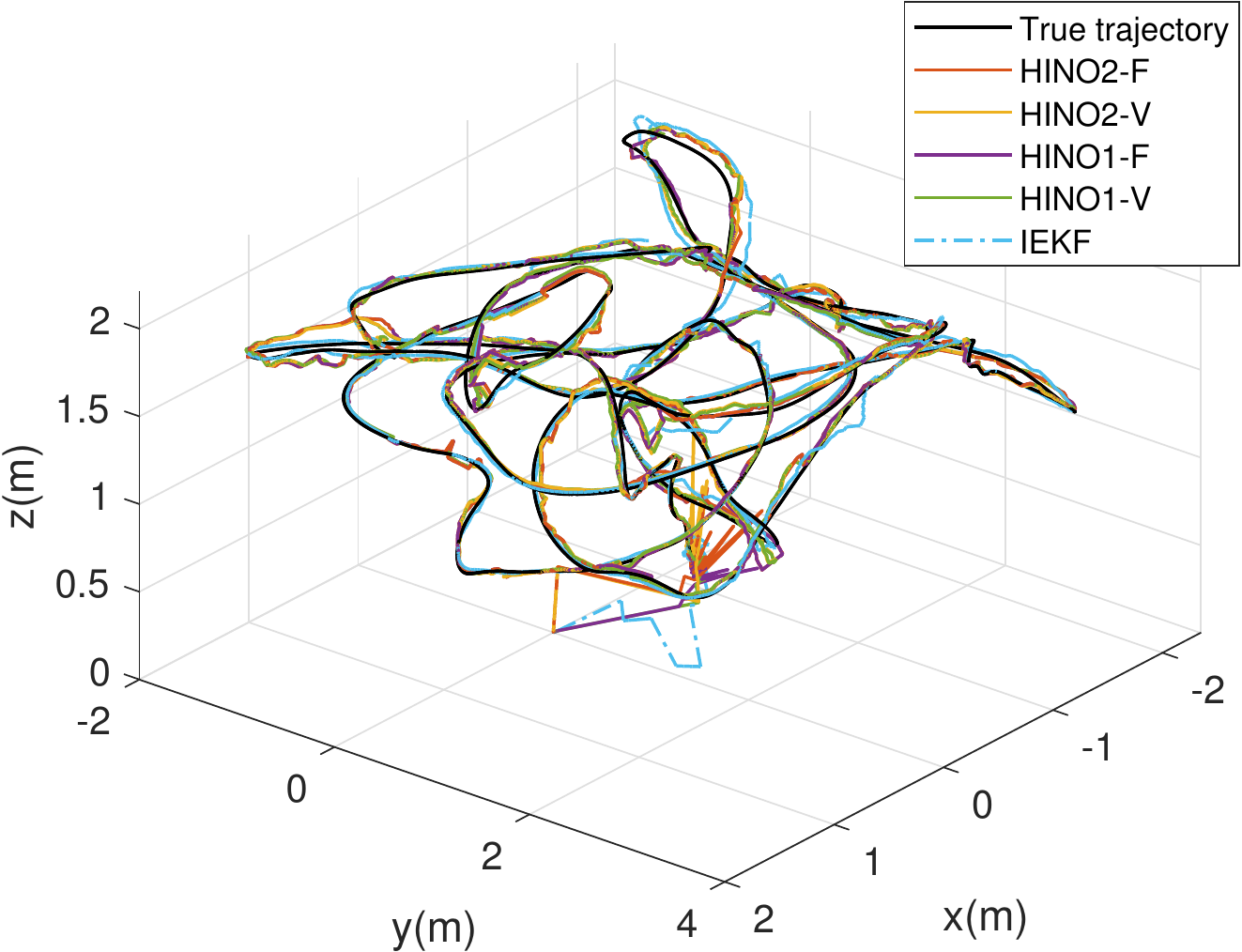}~
				\includegraphics[width=0.52\linewidth]{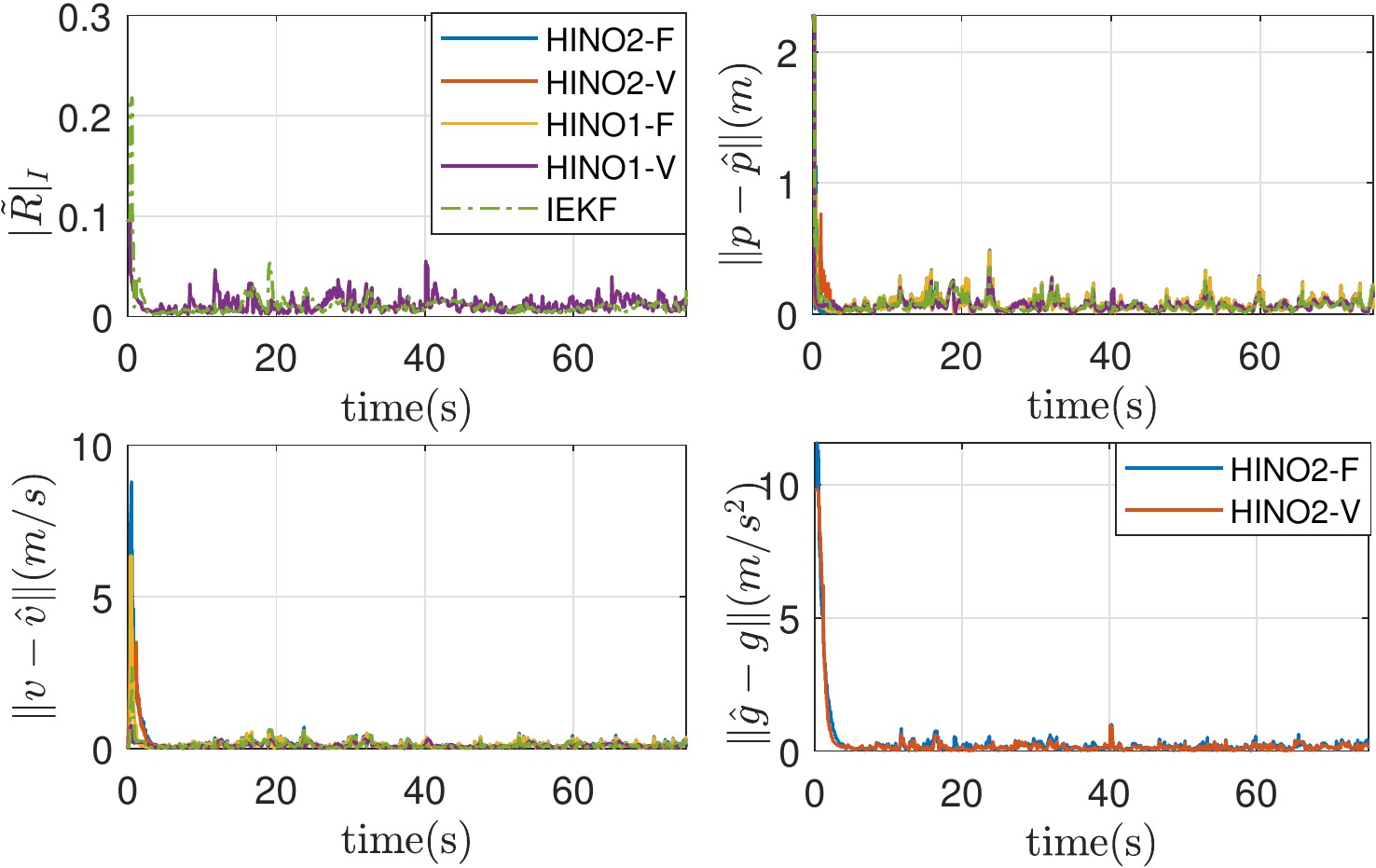}
			\end{minipage}\\[0.1cm]
		{\footnotesize  (a) Experimental results using dataset V1\_02\_medium }  \\[0.1cm]
		\begin{minipage}{0.92\linewidth}
			\includegraphics[width=0.40\linewidth]{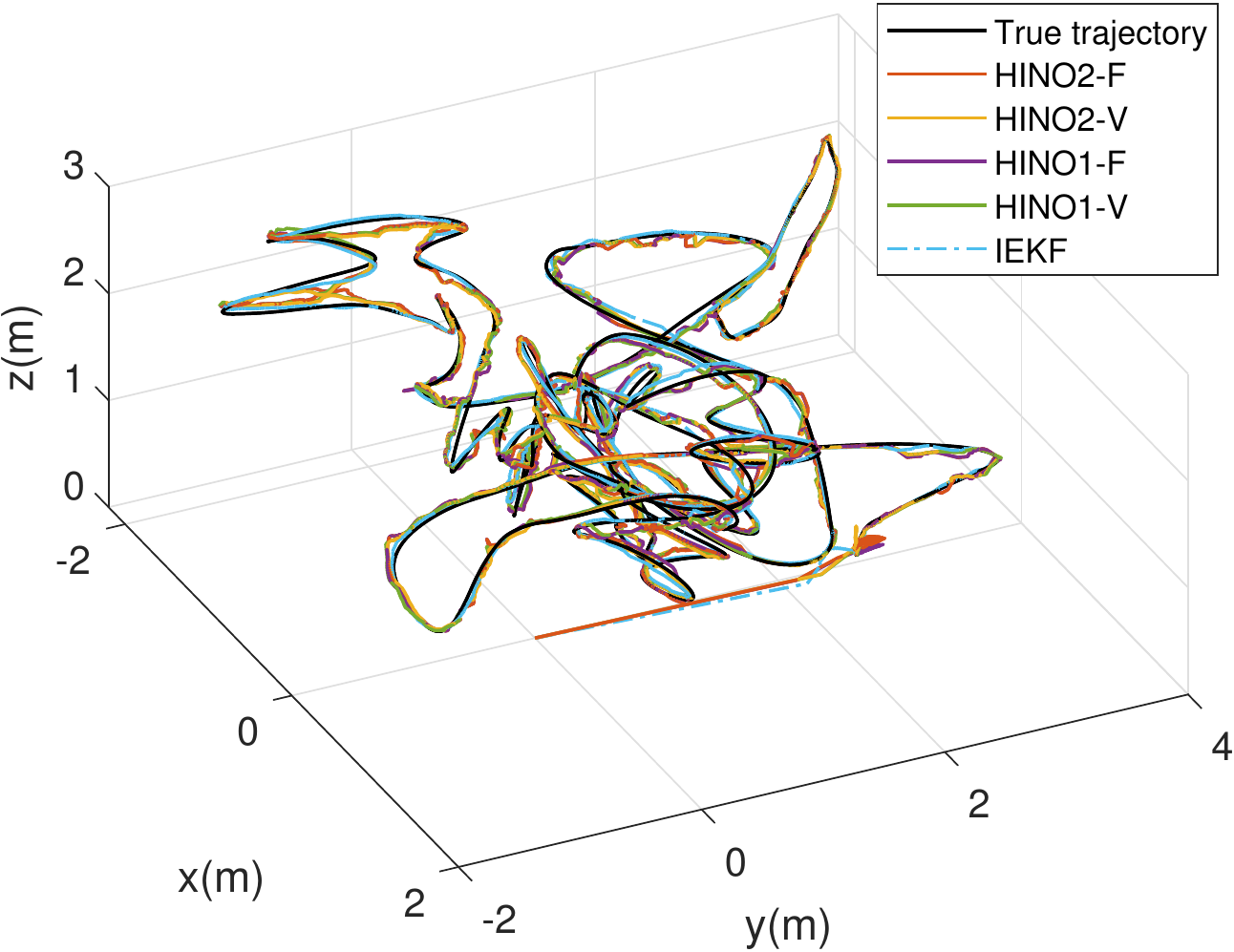}~
			\includegraphics[width=0.52\linewidth]{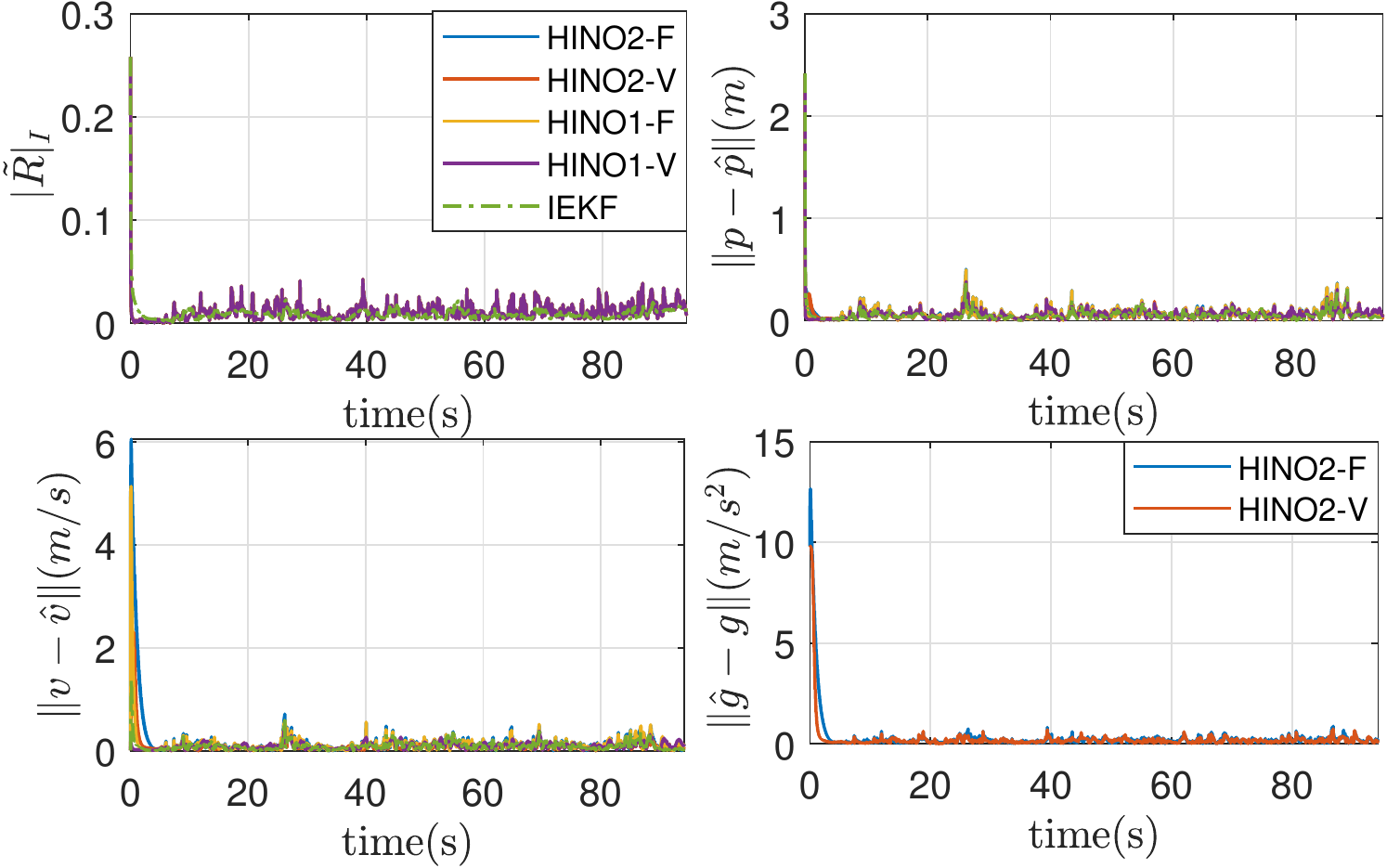}
		\end{minipage}\\[0.1cm]
		{\footnotesize  (b) Experimental results using dataset V1\_03\_difficult }  \\[0.1cm]
		\caption{The experimental results using the EuRoc dataset \cite{Burri25012016} with initial conditions: $\hat{R}(0)=\exp(0.1\pi u^\times)R_G$ with $ u\in \mathbb{S}^2$, $\hat{p}(0)=\hat{v}(0)= \eta(0)=\hat{g}(0) = 0_{3\times 1}$.
			The true (groundtruth) and estimated trajectories are shown in the left. The estimation errors of rotation, position, velocity and gravity vector are shown in the right.
		}
		\label{fig:experiment}
	\end{figure*}
	To achieve fast convergence for the attitude estimation, we choose $k_R=28/\|M\|_F$ with $M=\sum_{i=1}^N k_i (p_i-p_c)(p_i-p_c)\T$  and  $k_i=1/N$ for all $i=1,\hdots,N$  for each of the proposed observers. For the variable-gain observers, the matrices $V_t$ and $Q_t$ for the continuous-discrete Riccati equations are chosen same as the simulation with the covariance matrices of the measurement noise  $\text{Cov}(n_\omega) = 0.0024 I_3, \text{Cov}(n_a) = 0.0283 I_3$ and $ \text{Cov}(n_y^i) = 0.06 I_3$ for all $  i=1,\hdots N$. For the fixed-gain observers, to achieve similar performance as the variable-gain observers, the gains are picked after several  trials as $k_p = 0.85, k_v = 2.5$ and $k_{g} = 2.0$ with the condition (\ref{eqn:defXi_P}) satisfied.  Since the IMU measurements are not continuous although obtained at a high rate (200Hz), we  use the following numerical integration for the  estimated attitude: $\hat{R}_{k+1} = \hat{R}_k\exp(\omega_k^\times \Delta_T)$ for all $k\in \mathbb{N}_{>0}$ with the sampling period $\Delta_T= 0.005s$. The other state variables, namely $\hat{p}, \hat{v}, \eta$ and $\hat{g}$, are integrated using  Euler method with the sampling period  $\Delta_T$.
	The experimental results  are shown in   Fig. \ref{fig:experiment}. As one can see, the estimates provided by all the proposed hybrid observers and IEKF converge, after a few seconds, to the vicinity of the ground truth. The experimental results show that the performance of the proposed observers, which require less computational cost, are comparable to the performance of the IEKF.
	
	\section{Conclusion}
	Hybrid inertial navigation observers relying on continuous angular velocity and linear acceleration measurements, and intermittent landmark position measurements, have been proposed. Different versions have been developed depending on whether the observer gains are constant or time-varying and whether the gravity vector is known or not. While the fixed-gain observers require less computational cost than the variable-gain observers, the latter are easily tunable via the covariance of the measurements noise. All the proposed observers are endowed with strong exponential stability properties. This is, to the best of our knowledge, the first work dealing with inertial navigation observers design, using intermittent measurements, with strong stability guarantees. Simulation and experimental results, illustrating the performance of the proposed hybrid nonlinear observers, have been provided. 

	\begin{ack}                               
		The authors would like to thank Dr. Soulaimane Berkane for his interesting discussions with regards to Lemma \ref{lem:lemma-CDRE}.
	\end{ack}
	
	\section*{Appendix}	
	\appendix	
	\section{Solving the infinite-dimensional problem} \label{sec:solvingXoP}

	In this section, we provide a procedure motivated by \cite{sferlazza2019time} to solve the infinite-dimensional problem $\Xi_P(\tau)<0,$ for all $\tau\in [T_m,T_M]$. First, let  $(\lambda_i(\tau), v_i(\tau))\in \mathbb{R}\times \mathbb{S}^{n-1}$ be the $i$-th pair of eigenvalue and eigenvector of the matrix $\Xi_P(\tau)$ with $\lambda_i(\tau) = v_i(\tau)\T \Xi_P(\tau) v_i(\tau)$ being a continuous function of $\tau$. Using the facts   $\partial \Phi(\tau)/\partial \tau = e^{A\tau} A = Ae^{A\tau}$ and $ v_i(\tau)\T \left( \partial v_i(\tau)/\partial \tau \right)  = 0$ for all $\tau\geq 0$,  one has 	
	$
	\left|  \frac{\partial \lambda_i(\tau)}{\partial \tau} \right|  = \left| v_i\T   \frac{\partial\Xi_P(\tau)}{\partial \tau}   v_i \right|
	= v_i\T (A_g\T e^{A\T \tau}(A\T P + PA)e^{A\tau} A_g )v_i
	\leq 2 \|P\| \|A\| \|e^{A \tau}\|^2 \|A_g\|^2 : =\delta^*_P
	$.  Then, according to \cite[Lemma 4]{sferlazza2019time}, if there exist a constant $\mu>0$ and a scalar $\bar{\tau}\in [T_m,T_M]$ such that $\Xi_P(\bar{\tau})<-2\mu I$, then the maximum eigenvalue of $\Xi_P(\tau)$ cannot be greater than $-\mu$ as long as $\tau\in [\bar{\tau} - \mu/\delta^*_P  , \bar{\tau}+ \mu/ \delta^*_P ]$.
	
	The following procedure adapted from \cite[Algorithm 1]{sferlazza2019time} is presented to solve the infinite-dimensional problem  in a finite number of steps:
	\begin{itemize}
		\item[] \textbf{Step 1:} Obtain an exponential bound for $e^{A \tau}, \forall \tau\in [T_m, T_M]$ by finding a solution $\Pi=\Pi\T>0$ and $\beta >0$ satisfying
		\begin{equation}
		(-A + \beta I)\T \Pi + \Pi(-A + \beta I) >0 \label{eqn:solveII}.
		\end{equation}
		From \cite[Lemma 3]{sferlazza2019time}, one can show that $\|e^{A \tau}\| \leq   \gamma e^{\beta \tau}\leq   \gamma e^{\beta T_M}:= c_{A}, \forall \tau\in [T_m,T_M]$ with $  \gamma: = \sqrt{ {\lambda_M(\Pi)}/{\lambda_m(\Pi)}}$.
		
		\item[]  \textbf{Step 2:} Solve the finite-dimensional optimization problem with a constant $\mu>0$ and a discrete set $\mathcal{T}\subset [T_m, T_M]$(in the first step $\mathcal{T}= \{T_m, T_M\}$)
		\begin{align}
		(P^*,\bar{p}^*)  = & \arg \min_{P=P\T,\bar{p}} \bar{p},   \text{ subject to }  \nonumber \\
		& \Xi_P(\tau) < -2\mu I, \tau\in \mathcal{T} \nonumber \\
		& I\leq P \leq \bar{p} I  . \label{eqn:finiteoptimization}
		\end{align}	
		\item[] \textbf{Step 3:} Let $\delta^*_\mathcal{T} := \mu(2 \bar{p}^*  c_{A}^2 \|A\|\|I-CK\|^2)^{-1}$, and define a finite discrete set $\mathcal{T}_d=[T_m:  2\delta^*_\mathcal{T} :T_M]$. Check the eigenvalue condition $\Xi_{P^* }(\tau)<-\mu I, \forall  \tau\in \mathcal{T}_d$. 
		If this step is successful, the solution $P^*$ from Step 2 is a solution of the infinite-dimensional problem $\Xi_P(\tau)<0$ for all $\tau\in [T_m,T_M]$ and then the algorithm stops. Otherwise, add the worst-case value $\tau = \arg \max_{\tau\in \mathcal{T}_d} (\lambda_M(\Xi_{P^*}(\tau)))$ to the discrete set $\mathcal{T}$ and restart from Step 2 again.
	\end{itemize}
	Note that if the infinite-dimensional problem $\Xi_P(\tau)<0, \forall \tau\in [T_m,T_M]$ is not feasible, one has to redesign the gain matrix $K$  such that this optimization problem is feasible. The finite-dimensional optimization problem (\ref{eqn:solveII}) and (\ref{eqn:finiteoptimization})  can be solved by a
	convex optimization solver like CVX   \cite{grant2009cvx}.
	
	\section{Proof of Theorem \ref{theo:1}} \label{sec:theo1}
	Before proceeding with the proof of Theorem 1, we provide the
	following useful properties: for any  solution of $\dot{ {R}} =  {R}(\omega)^\times$ with $R(0)\in SO(3)$ and $\omega \in  \mathbb{R}^3$, one has
	\begin{align}
	&4\lambda_{m}^{\bar{M}}| {R}|_I^2 \leq  \tr((I- {R})M) \leq 4  \lambda_{M}^{\bar{M}}| {R}|_I^2 \label{eqn:property1}\\
	&\|\psi(M {R})\|^2 =  \alpha(M, {R})\tr((I- {R})W') \label{eqn:property2}\\
	& \dot{\psi}(M {R}) =  E(M {R}) \omega  \label{eqn:property3}\\
	& \|E(M {R})\|_F \leq \|\bar{M}\|_F\label{eqn:property4}
	\end{align}
	where $M=M\T\geq 0$, $E(M {R}) := \frac{1}{2}(\tr(M\tilde{R})I - \tilde{R}\T M)$, $\bar{M} = \frac{1}{2}(\tr(M)I_3 - M)$, $W' := \tr(\bar{M}^2)I - 2\bar{M}^2$, and $\alpha(M, {R}):=1-| {R}|_I^2 \cos^2\langle u,\bar{M}u\rangle  \leq 1 $ with $\langle~,~\rangle$ denoting the angle between two vectors and $u$ being the axis of
	rotation of $ {R}$.  The proof of these properties can be found in \cite[Lemma 1 and Lemma 2]{berkane2017hybrid} and the references therein.

	First, we are going to show that $|\tilde{R}(t)|_I < 1$ for all $ t \geq 0$. This step guarantees that the innovation term $\sigma_R = k_R \psi(M\tilde{R})$ vanishes only at $\tilde{R}=I_3$ excluding the undesired critical points  $\tilde{R} = \mathcal{R}_a(\pi, v)$ with $v\in \mathcal{E}(M)$.
	Consider the following real-valued function on $\mho \times [0,T_M]$:
	\begin{align}
	\mathcal{W} =  \frac{1}{2} \tr((I-\tilde{R})M) -  \tau\eta\T\psi(M\tilde{R}) +   \mu e^{\tau} \eta\T \eta \label{eqn:defW}
	\end{align}
	with some $\mu>0$ .  Let $\|\zeta\|^2 := |\tilde{R}|_I^2 +  \|\eta\|^2$. Using (\ref{eqn:property1}), (\ref{eqn:property2}), one has $\|\psi(M\tilde{R})\|^2 \leq   \tr((I-\tilde{R})W') \leq 4  \lambda_{M}^{\bar{W}}  |\tilde{R}|_I^2  $ with $\bar{W}: = \frac{1}{2}(\tr(W')I_3-W')$. For all $\tau \geq 0$, one obtains the following inequalities:
	\begin{equation}
	\lambda_m^{P_1 } \|\zeta\|^2 \leq \mathcal{W} \leq \lambda_M^{P_2 }\|\zeta\| ^2  \label{eqn:boundV_1}
	\end{equation}
	where   the  matrices $P_1 $ and $ P_2 $ are given by
	\begin{align*} 
	P_1 = \begin{bmatrix}
	2\lambda_{m}^{\bar{M}}  & - c_{\bar{W}} T_M  \\
	- c_{\bar{W}} T_M  &  \mu
	\end{bmatrix},   
	P_2   = \begin{bmatrix}
	2\lambda_{M}^{\bar{M}}  & c_{\bar{W}} T_M \\
	c_{\bar{W}} T_M &  \mu e^{T_M}
	\end{bmatrix}
	\end{align*}
	with $c_{\bar{W}} = \sqrt{\lambda_{M}^{\bar{W}} }$.
	To guarantee that   $P_1$ and $P_2$ are positive definite, it is sufficient to choose   $\mu>   \frac{1}{2}\lambda_{M}^{\bar{W}} T_M^2 /\lambda_{m}^{\bar{M}} $.  Then, the time-derivative of $\mathcal{W}$ along the flows of (\ref{eqn:closed-loop1}) is given by
	\begin{align}
	\dot{\mathcal{W}} & \leq   \frac{1}{2} \tr(M\tilde{R}\eta^\times ) + \eta\T\psi(M\tilde{R}) +  \tau \eta\T E(M\tilde{R}) \eta - \mu   \|\eta\|^2 \nonumber \\
	& \leq    \left(  T_M \|\bar{M}\|_F - \mu  \right)  \eta\T\eta   \label{eqn:dotW_1}
	\end{align}
	where we made use of  (\ref{eqn:property3}) and (\ref{eqn:property4}), and the facts $e^{\tau} \geq 1$ for all $ \tau \geq 0$ and $\tr(M\tilde{R}\eta^\times )= -2 \eta \T \psi(M\tilde{R})$. Choosing
	$
	\mu >  \max  \{ \frac{1}{2}\lambda_{M}^{\bar{W}} T_M^2 /\lambda_{m}^{\bar{M}}, T_M \|\bar{M}\|_F\} ,
	$
	one obtains $\dot{\mathcal{W}} \leq  0$,
	which implies that $\dot{\mathcal{W}}$ is negative semi-definite and  $\mathcal{W}$ is non-increasing in the flows.

	Thanks to the decreasing timer  in (\ref{eqn:tau}), one has $\tau = 0$ at each jump. Let $\mathcal{W}^+$ be the value of $\mathcal{W}$ after each jump. Then,  one can show that
	\begin{align}
	\mathcal{W}^+ - \mathcal{W} & =  \frac{1}{2}\tr((I-\tilde{R}^+)M)  - \frac{1}{2}  \tr((I-\tilde{R})M)  \nonumber\\
	&\quad   -  \nu(\eta^+)\T\psi(M\tilde{R}^+) +    \tau  \eta \T\psi(M\tilde{R} ) \nonumber  \\
	& \quad +    \mu e^{\nu}  \|\eta^+\|^2 -  \mu  \|\eta\|^2 \nonumber\\
	& \leq  -  k_R \left(   T_m    - \mu  e^{T_M} k_R   \right) \| \psi(M\tilde{R})\|^2  -  \mu \|\eta\|^2 \nonumber \\
	& \leq   - k_R  \varrho^* \|\psi(M\tilde{R})\|^2-   \mu \|\eta\|^2   \label{eqn:Wdecrease}
	\end{align}
	where $\nu := \tau^+ \in [T_m,T_M]$, $\varrho^*: =  T_m    - \mu  e^{T_M} k_R $, and  we made use of the fact that $\tilde{R}^+= \tilde{R}$ and $ \eta^+ = k_R \psi(M\tilde{R})$. Choosing
	$ k_R <   k_R^* :=\frac{1}{\mu} T_m   e^{-T_M} $ such that $\varrho^* =  T_m    - \mu  e^{T_M} k_R >0$,
	one can further show that $\mathcal{W}^+ - \mathcal{W} \leq 0$,
	which implies that $\mathcal{W}$ is non-increasing after each jump.
	From (\ref{eqn:boundV_1}), (\ref{eqn:dotW_1}) and (\ref{eqn:Wdecrease}), one can show that $\lambda_m^{P_1} |\tilde{R}|_I^2  \leq \lambda_m^{P_1} \|\zeta\|^2  \leq \mathcal{W}(t) \leq \mathcal{W}(0)$ for all $t \geq 0$. The minimum eigenvalue of $P_1$ is explicitly given by $	0<\lambda_m^{P_1}
	= \lambda_m^{\bar{M}} + \frac{ \mu}{2} - \frac{1}{2}( (2 \lambda_m^{\bar{M}} - \mu)^2 + 4 \lambda_{M}^{\bar{W}} T_M^2)^{\frac{1}{2}}< 2\lambda_m^{\bar{M}}$ with $\mu >  \max  \{ \frac{1}{2}\lambda_{M}^{\bar{W}} T_M^2 /\lambda_{m}^{\bar{M}}, T_M \|\bar{M}\|_F\} $.
	It is easy to verify that $\lambda_m^{P_1}$ is a continuous monotonically increasing function of $\mu>0$ and  $\lim_{\mu \rightarrow +\infty} \lambda_m^{P_1}=  2 \lambda_m^{\bar{M}}$.
	Hence, given a constant $0< \epsilon^* <   1-\epsilon^2 $, it is always possible to find a constant $\mu_{\epsilon}^* >0$ (depending on $\epsilon$ and $\epsilon^*$) such that, for any $ \mu  >    \mu_{\epsilon}^* $,
	one has $\lambda_m^{P_1}  \geq 2 (\epsilon^2 + \epsilon^*) \lambda_m^{\bar{M}} = 2 (\epsilon^2 + \epsilon^*) \varsigma_{M} \lambda_M^{\bar{M}}$.   From (\ref{eqn:boundV_1}), substituting  $|\tilde{R}(0)|_I^2 \leq \epsilon^2 \varsigma_{M} < 1$ and $\|\eta(0)\|\leq \eta^*$, one can show that
	$
	\mathcal{W}(0) \leq   2 \lambda_M^{\bar{M}} |\tilde{R}(0)|_I^2  +  2T_M c_{\bar{W}} \|\eta(0)\|   +   \mu e^{T_M} \|\eta(0)\|^2
	\leq  2 \lambda_M^{\bar{M}}  ( \epsilon^2 \varsigma_{M}   +   h(\eta^*)  )
	$, 
	where $h(\eta^*): = (T_M  c_{\bar{W}} \eta^*  +  \frac{1}{2} \mu e^{T_M} (\eta^*)^2)/{ \lambda_M^{\bar{M}}}$. Since  $h(\eta^*)$ is a continuous monotonically increasing function of $\eta^*\geq 0$ and $h(0)=0$, there exists a constant $\eta^*: = h^{-1}(c_\eta \epsilon^* \varsigma_{M})>0$  with some constant $0<  c_\eta<1$ such that for all $\|\eta(0)\|\leq \eta^*$ one has $\mathcal{W}(0) \leq 2 \lambda_M^{\bar{M}} (\epsilon^2 +   c_\eta \epsilon^*) \varsigma_{M}$. Hence, one can conclude that 
	\begin{equation} 
	|\tilde{R}|_I^2   \leq \frac{1}{\lambda_m^{P_1} }  \mathcal{W}(0)  \leq  \frac{\epsilon^2 +   c_\eta \epsilon^* }{   \epsilon^2 + \epsilon^*}    <   1   \label{eqn:Rless1}
	\end{equation}
	for all $t \geq 0$ with $k_R <   k_R^*  =\frac{1}{\mu} T_m   e^{-T_M}$   and $\mu >  \max  \{\mu_{\epsilon}^*, \frac{1}{2}\lambda_{M}^{\bar{W}} T_M^2 /\lambda_{m}^{\bar{M}}, T_M \|\bar{M}\|_F\}$. 
	
	Next, let us show the exponential stability of the set $\mathcal{A}$.
	From (\ref{eqn:Rless1}), one obtains $|\tilde{R}|_I^2 \leq 1- \frac{ (1- c_\eta) \epsilon^* }{   \epsilon^2 + \epsilon^*}  $ for all $t \geq 0$. Then, applying (\ref{eqn:property1}) and (\ref{eqn:property2}) one has $\|\psi(M\tilde{R})\|^2  \geq 4 (1-|\tilde{R}|_I^2)\lambda_{m}^{\bar{W}}|\tilde{R}|_I^2 \geq  4\frac{ (1- c_\eta) \epsilon^* }{   \epsilon^2 + \epsilon^*} \lambda_{m}^{\bar{W}} |\tilde{R}|_I^2$. From (\ref{eqn:boundV_1}) and (\ref{eqn:Wdecrease}),   it follows that 
	\begin{align}
	\mathcal{W}^+ &\leq \mathcal{W}-  4k_R \varrho^* \lambda_m^{\bar{W}} \frac{ (1- c_\eta) \epsilon^* }{   \epsilon^2 + \epsilon^*}   |\tilde{R}|_I^2 -  \mu \|\eta\|^2 \nonumber \\
	&\leq \left(1 - \frac{c_1}{\lambda_M^{P_2}} \right)\mathcal{W} =  e^{-c_2} \mathcal{W} \label{eqn:V_1+}
	\end{align}
	with $c_1:=\min\{ 4k_R \varrho^* \lambda_m^{\bar{W}}\frac{ (1- c_\eta) \epsilon^* }{   \epsilon^2 + \epsilon^*}   ,  \mu \} $ and $c_2: =  -\ln (1 -  {c_1}/{\lambda_M^{P_2}}  )$.   Since $\mathcal{W}^+  \geq 0$ and ${c_1}/{\lambda_M^{P_2}}>0$, one can show that $0<1 -  {c_1}/{\lambda_M^{P_2}} <1$ and then $c_2>0$ is well-defined. Consider the following real-valued function on $SO(3)\times \mathbb{R}^3 \times[0, T_M]$:
	\begin{equation}
	\mathcal{V}_1(\tilde{R},\eta,\tau) = e^{\lambda_{1}^F\tau} \mathcal{W} \label{eqn:V1}
	\end{equation}
	with $0< \lambda_{1}^F< c_2/T_M$.
	In view of (\ref{eqn:boundV_1}), (\ref{eqn:dotW_1}) and (\ref{eqn:V_1+}), one can show that
	\begin{align}
	\underline{\alpha}_1\|\zeta\|^2  &\leq \mathcal{V}_1(\tilde{R},\eta,\tau) \leq  \bar{\alpha}_1   \|\zeta\|^2  \label{eqn:V1bound}\\
	\dot{\mathcal{V}}_1(\tilde{R},\eta,\tau) &\leq -\lambda_{1}^F\mathcal{V}_1(\tilde{R},\eta,\tau), \quad \forall x_1\in \mathcal{F}_1 \label{eqn:dotV1}\\
	\mathcal{V}_1(\tilde{R}^+,\eta^+,\tau^+) &\leq e^{-\lambda_{1}^J} \mathcal{V}_1(\tilde{R},\eta,\tau),\quad \forall x_1\in \mathcal{J}_1 \label{eqn:V1+}
	\end{align}
	where $ \underline{\alpha}_1:=  \lambda_m^{P_1}  , \bar{\alpha}_1:=e^{(\lambda_{1}^F T_M)} \lambda_M^{P_2}, \lambda^J_{1}:= c_2-\lambda_{1}^F T_M>0$, and we made use of the facts: $\mathcal{V}_1^+ = e^{ \lambda_{1}^F \tau^+} \mathcal{W}^+ \leq  e^{( \lambda_{1}^F T_M -c_2)} \mathcal{W} $ and $\mathcal{W} \leq \mathcal{V}_1$. From (\ref{eqn:V1bound})-(\ref{eqn:V1+}), one can conclude that $\mathcal{V}_1$ is exponentially decreasing in both flow and jump sets, which further implies that the estimation error $\zeta$ converges exponentially to zero.

	On the other hand, let us consider the following real-valued function on $\mathbb{R}^6 \times [0,T_M]$:
	\begin{equation}
	\mathcal{V}_2(\mathsf{x},\tau) = e^{\lambda^F_2 \tau}\mathsf{x}\T \Phi\T(\tau)P \Phi(\tau) \mathsf{x} \label{eqn:V2},
	\end{equation}
	with some $\lambda^F_2>0$.  One can easily verify that there exist two positive constants $\underline{\alpha}_2,\bar{\alpha}_2$ such that
	\begin{equation}
	\underline{\alpha}_2 \|\mathsf{x}\|^2 \leq \mathcal{V}_2(\mathsf{x},\tau)  \leq  \bar{\alpha}_2 \|\mathsf{x}\|^2 \label{eqn:V2-bound}
	\end{equation}
	with $
	\underline{\alpha}_2   :=  \min_{\tau \in [0, T_M]} e^{\lambda^F_2 \tau} \lambda_{m}^{\left(   \Phi\T(\tau)P \Phi(\tau)\right)  }  $ and $
	\bar{\alpha}_2   :=  \max_{\tau \in [0, T_M]} e^{\lambda^F_2 \tau} \lambda_{M}^{ \left( \Phi\T(\tau)P \Phi(\tau) \right) }.
	$
	The time-derivative of $\mathcal{V}_2$ along the flows of (\ref{eqn:closed-loop1}) is given by
	\begin{align}
	\dot{\mathcal{V}}_2 & = - \lambda_2^F \mathcal{V}_2 + e^{\lambda^F_2 \tau}\mathsf{x}\T \Phi\T(\tau)(A \T P + P A   ) \Phi(\tau) \mathsf{x} \nonumber \\
	&\quad + 2 e^{\lambda^F_2 \tau}\mathsf{x}\T \Phi\T(\tau) P \Phi(\tau) \delta_g  \nonumber \\
	& \quad  + e^{\lambda^F_2 \tau}\mathsf{x}\T (\dot{\Phi}\T(\tau) P\Phi(\tau) + \Phi\T(\tau) P \dot{\Phi}(\tau))  \mathsf{x} \nonumber \\
	& \leq  - \lambda_2^F \mathcal{V}_2 + 2 \bar{\alpha}_2 \|\mathsf{x}\| \|\delta_g \| \nonumber \\
	& \leq  - \lambda_2^F \mathcal{V}_2 +  \beta \|\mathsf{x}\| \|\zeta\|, \quad \forall x_1\in \mathcal{F}_1  \label{eqn:dotV2}
	\end{align}
	where $\beta:= 4\sqrt{2} c_g  \bar{\alpha}_2$, and we made use of the facts $\|\delta_g \|\leq c_g  \|I- \tilde{R}\|_F =2\sqrt{2} c_g  |\tilde{R}|_I \leq 2\sqrt{2} c_g  \|\zeta\|$, $A\Phi(\tau) = \Phi(\tau) A$, and $\dot{\Phi}(\tau) = -A\Phi(\tau)$. Since $\Xi_P(\tau)<0, \forall \tau\in[T_m,T_M]$, there exists a (small enough) positive scalar $c_q < p_M:=\lambda_{M}^P$ such that $\Xi_P(\tau) \leq -c_q I, \forall \tau\in[T_m,T_M]$. Let $\mathcal{V}_2^+:= \mathcal{V}_2(\mathsf{x}^+,\tau^+)$ and $\nu:=\tau^+\in [T_m,T_M]$. Then, for each jump at $\tau=0$, one has
	\begin{align}
	\mathcal{V}_2^+ 
	&= 	e^{\lambda^F_2 \nu} \mathsf{x}\T ( (I-KC)\T\Phi(\nu)\T P \Phi(\nu) (I-KC)  ) \mathsf{x} \nonumber \\
	&= 	e^{\lambda^F_2 \nu} \mathsf{x}\T   \Xi_P(\nu) \mathsf{x} +  e^{\lambda^F_2 \nu}  \mathsf{x}\T  P \mathsf{x}\nonumber \\
	&\leq - e^{\lambda^F_2 \nu} c_q \mathsf{x}\T  \mathsf{x} +  e^{\lambda^F_2 \nu}   p_M \mathsf{x}\T  \mathsf{x} \nonumber \\
	&\leq \left( p_M-c_q  \right) e^{\lambda^F_2 T_M} \|\mathsf{x}\|^2 \nonumber \\
	& \leq e^{-\lambda_2^J}\mathcal{V}_2(\mathsf{x} ), \quad \forall x_1\in \mathcal{J}_1 \label{eqn:V2+}
	\end{align}
	where  $ \lambda_2^J: =  -\ln (\frac{p_M-c_q}{\bar{\alpha}_2} e^{\lambda^F_2 T_M} )$. Pick $\lambda_2^F  < \frac{1}{T_M}\ln  (\frac{\bar{\alpha}_2}{p_M-c_q}  )$ such that $  \frac{p_M-c_q}{\bar{\alpha}_2} e^{\lambda^F_2 T_M} <1$ and $\lambda_2^J>0$.
	
	Now, we are going to show the exponential stability of the set $ {\mathcal{A}}$ for the overall hybrid closed-loop system $\mathcal{H}_1$ in (\ref{eqn:closed-loop1}). Let  $|x_1|_{\mathcal{A}} \geq 0$ be the distance of $x_1$ with respect to the set $\mathcal{A}$ such that $|x_1|_{\mathcal{A}}^2:=\inf_{(\bar{R},\bar{\eta},\bar{\mathsf{x}},\bar{\tau})\in \mathcal{A}} (  |\tilde{R}\bar{R}\T|_I^2 + \|\eta-\bar{\eta}\|^2 + \|\mathsf{x}- \bar{\mathsf{x}}\|^2 + \|\tau-\bar{\tau}\|^2) = |\tilde{R}|_I^2 + \|\eta\|^2 + \|\mathsf{x}\|^2 = \|\zeta\|^2 +   \|\mathsf{x}\|^2 $. Consider the Lyapunov function candidate $\mathcal{V}(x_1)= \varepsilon \mathcal{V}_1 (\tilde{R},\eta,\tau) + \mathcal{V}_2(\mathsf{x},\tau)$ with some $\varepsilon>0$. From (\ref{eqn:V1bound}) and (\ref{eqn:V2-bound}), one can show that
	\begin{align}
	& \underline{\alpha} |x_1|_{ {\mathcal{A}}}^2 \leq \mathcal{V}(x_1) \leq \bar{\alpha} |x_1|_{ {\mathcal{A}}}^2 \label{eqn:Vbound}
	\end{align}
	where $\underline{\alpha}:=\min\{\varepsilon\underline{\alpha}_1, \underline{\alpha}_2\}$ and $ \bar{\alpha}:=\max\{\varepsilon\bar{\alpha}_1,\bar{\alpha}_2\}$.
	In view of (\ref{eqn:V1})-(\ref{eqn:dotV1}) and  (\ref{eqn:V2-bound})-(\ref{eqn:dotV2}),  one has
	\begin{align}
	\dot{\mathcal{V}} & = -\varepsilon\lambda_{1}^F\mathcal{V}_1 - \lambda_2^F \mathcal{V}_2 +  \beta \|\mathsf{x}\| \|\zeta\| \nonumber \\
	& \leq -\varepsilon \lambda_{1}^F \mathcal{V}_1 - \lambda_{2}^F \mathcal{V}_2 + \frac{\varepsilon \lambda_{1}^F \underline{\alpha}_1}{2}  \|\zeta\|^2 + \frac{\beta^2}{2\varepsilon \lambda_{1}^F \underline{\alpha}_1} \|\mathsf{x}\|^2 \nonumber \\ 
	& \leq -\frac{\varepsilon}{2} \lambda_{1}^F \mathcal{V}_1 - \left(  \lambda_{2}^F  - \frac{\beta^2}{2\varepsilon \lambda_{1}^F \underline{\alpha}_1 \underline{\alpha}_2}\right)  \mathcal{V}_2 \nonumber \\
	&= -\lambda_F \mathcal{V}, \quad \forall x_1\in \mathcal{F}_1 \label{eqn:dotV}
	\end{align}
	with $	\varepsilon > \frac{\beta^2} {2\lambda_{1}^F \lambda_{2}^F \underline{\alpha}_1 \underline{\alpha}_2}$ such that $\lambda_F :=   \min\{ \frac{  \lambda_{1}^F  }{2} ,  ( \lambda_{2}^F -  \frac{\beta^2}{2\varepsilon \lambda_{1}^F \underline{\alpha}_1 \underline{\alpha}_2} ) \} >0$, and we have made use of the inequality $\beta \|\mathsf{x}\| \|\zeta\| \leq  \frac{\varepsilon \lambda_{1}^F \underline{\alpha}_1}{2}  \|\zeta\|^2 + \frac{\beta^2}{2\varepsilon \lambda_{1}^F \underline{\alpha}_1} \|\mathsf{x}\|^2 $. From (\ref{eqn:V1+}) and (\ref{eqn:V2+}), one obtains
	\begin{equation}
	\mathcal{V}(x_1^+) 
	\leq e^{- \lambda_J} \mathcal{V} (x_1) , \quad \forall x_1\in \mathcal{J}_1 \label{eqn:V+}
	\end{equation}
	where $\lambda_J:=  \min\{\lambda_1^J,\lambda_2^J\}$. Let $\lambda_c: = \min\{\lambda_F,\lambda_J\}$. In view of  (\ref{eqn:dotV}) and (\ref{eqn:V+}), one has $\mathcal{V}(x_1(t,j))\leq e^{-\lambda_c(t+j)} \mathcal{V}(x_1(0,0))$ for all $(t,j)\in \dom x_1$. Since $T_m$ is strictly positive by Assumption \ref{assum:1}, it follows that every maximal solution to the hybrid system  $\mathcal{H}_1$ is complete  as $\dom x_1$ is unbounded \ie,  $t+j \to \infty$. Applying (\ref{eqn:Vbound}),  one can conclude that  for all $(t,j)\in \dom x_1$
	$
	|x_1(t,j)|_{ {\mathcal{A}}} \leq \sqrt{\frac{ \bar{\alpha}}{\underline{\alpha}}}  e^{-\frac{1}{2}\lambda_c(t+j)}|x_1(0,0)|_{ {\mathcal{A}}},  
	$
	which shows that the set $ {\mathcal{A}}$ is exponentially stable. This completes the proof.

	\section{Proof of Theorem \ref{theo:2}} \label{sec:theo2}
	
	Following the same steps as in the first part of the proof of Theorem \ref{theo:1}, one obtains $|\tilde{R}|_I<1, \forall t \geq 0$. Using the real-valued function $\mathcal{V}_1$ defined in (\ref{eqn:V1}), one obtains inequalities (\ref{eqn:V1bound})-(\ref{eqn:V1+}) and it follows that the estimation error $\zeta$ converges exponentially to zero. 
	Since $\omega(t)$ is continuous and bounded, it is clear that $A_t$ and its associated transition matrix $\Phi_s^t$ are continuous and bounded. As shown in the proof of \cite[Lemma 3]{wang2020hybrid}, $\Phi_s^t$ can be written as $\Phi^t_s = T(t)\bar{\Phi}^t_s T^{-1}(t)$ where $T = \text{blkdiag}(R,R)$ and $\bar{\Phi}^t_s=\exp(A(t-s))$ denoting the state transition matrix associated to $A$ in \eqref{eqn:AKC1}. Since $T(t) T^{-1}(t) = I$ and $\exp(At) \geq  I$ for all $t\geq 0$, one verifies that $(\Phi^{s}_{t_j})\T \Phi^{s}_{t_j} \geq I, \forall s\geq t_j$, which shows that \eqref{eqn:CDRE_condition1} is satisfied. Moreover, using similar steps as in the proof of  \cite[Lemma 3]{wang2020hybrid}, one obtains $W_o(t_j,t_{j+\Gamma})=\sum_{i=j}^{j+\Gamma}  (\Phi^{t_i}_{t_{j}})\T C\T  C \Phi^{t_i}_{t_{j}} = T(t_j) (\sum_{i=j}^{j+\Gamma}  (\bar{\Phi}^{t_i}_{t_{j}})\T C\T  C \bar{\Phi}^{t_i}_{t_{j}}) T^{-1}(t_j)$. One can also show that the matrix $[C\T,(C\bar{\Phi}^{t_i}_{t_{j}})\T]\T$ has full rank for all $t_i> t_j\geq 0$. Therefore, one concludes that there exist some positive constants ${\mu}_q \in \mathbb{R},\Gamma \in \mathbb{N}_{>0}$ such that $W_o(t_j,t_{j+\Gamma}) \geq \mu_q I$ for all $t_j\geq 0$. This, with the fact that  $Q_t$ is bounded and uniformly positive definite, implies the existence of the lower bound in \eqref{eqn:CDRE_condition3}. Moreover, since $\Phi_s^t, V_t$ are continuous and bounded, $Q_t$ is bounded, and $Q_t,V_t$ are uniformly positive definite, one can easily show the existence of the upper bounds in \eqref{eqn:CDRE_condition2} and \eqref{eqn:CDRE_condition3}, and the lower bound in  \eqref{eqn:CDRE_condition2}.
	Consequently, all the conditions in  Lemma \ref{lem:lemma-CDRE} are satisfied, which guarantees that a solution $P$ for (\ref{eqn:CDRE-C})-(\ref{eqn:CDRE-D}) exists, and there exist constants $0< p_m \leq p_M < \infty $  such that $p_m I \leq P \leq p_M I, \forall t \geq 0$.
	
	Now, consider the real-valued function $\mathcal{V}_2(\mathsf{x},\tau) = e^{-  \gamma \tau} \mathsf{x}\T P^{-1} \mathsf{x}$ with some $  \gamma>0$, whose upper and lower bounds are given by
	\begin{equation}
	\underline{\alpha}_2'  \|\mathsf{x}\|^2  \leq \mathcal{V}_2(\mathsf{x},\tau)  \leq \bar{\alpha}_2' \|\mathsf{x}\|^2 \label{eqn:V2-bound_2}
	\end{equation}
	where $\underline{\alpha}_2':=\frac{1}{p_M } e^{-  \gamma T_M} $ and $\bar{\alpha}_2': = \frac{1}{p_m } $. The time-derivative of $\mathcal{V}_2$ along the flows of (\ref{eqn:closed-loop2}) is given as
	\begin{align}
	\dot{\mathcal{V} }_2 &=    \gamma\mathcal{V}_2 +  e^{-  \gamma \tau} \mathsf{x}\T (A_t\T P^{-1} + P^{-1}A_t + \dot{P}^{-1}) \mathsf{x} \nonumber \\
	&\quad + 2 e^{-  \gamma \tau} \mathsf{x}\T  P^{-1} \bar{\delta}_g \nonumber \\
	& \leq \left(  \gamma - e^{-  \gamma T_M} \frac{v_m p_m}{p_M^2} \right)\mathcal{V}_2  + \frac{2}{p_m}  \|\mathsf{x}\|  \|\bar{\delta}_g \| \nonumber \\
	&= -\lambda_2^F \mathcal{V}_2  +  \beta \|\mathsf{x}\|  \|\zeta\| , \quad \forall x_2\in \mathcal{F}_2 \label{eqn:dotV2_4}
	\end{align}
	where $\lambda^F_2 := -  \gamma + e^{-  \gamma T_M}  \frac{v_m p_m}{p_M^2}> 0$ with $  \gamma$ small enough, $v_m: = \inf_{t \geq 0} \lambda_m^V$,   $\beta:= 4 \sqrt{2}  c_g /p_m $, and we made use of the fact $\|\bar{\delta}_g \|\leq c_g  \|I- \tilde{R}\|_F \leq 2\sqrt{2} c_g  \|\zeta\|$. Let $\varpi(  \gamma)= -   \gamma + e^{-  \gamma T_M} \frac{v_m p_m  }{p_M^2}$. Note that $\varpi(0)=\frac{v_m p_m}{p_M^2}>0$, $\varpi(  \gamma)$  decreases as $  \gamma$ increases, and the solution for $\varpi(  \gamma)=0$ is given by $  \gamma =\frac{1}{T_M}W(\frac{v_m p_m T_M }{p_M^2})$ with $W(\cdot)$ denoting the Lambert $W$ function. Hence, one can show that $\lambda^F_2 > 0$ for any $0<  \gamma < \frac{1}{T_M}W(\frac{v_m p_m T_M }{p_M^2})$. 
	Since the solution of $P$ is well defined for all $(t,j)\in \dom x_2$ and $(P^+)^{-1} = P^{-1}+ C\T Q_t C$, one verifies that $I-KC$ is full rank and $(P^+)^{-1}$ can be rewritten as $ (P^+)^{-1} = P^{-1}(I-KC)^{-1}$. Let $\mathcal{V}_2^+:= \mathcal{V}_2(\mathsf{x}^+,\tau^+)$ and $\nu:=\tau^+\in [T_m,T_M]$. Then, for each jump at $\tau =0$, one has
	\begin{align}
	\mathcal{V}_2^+
	& =   e^{-  \gamma \nu }  \mathsf{x} \T (I-KC)\T (P^{+})^{-1} (I-KC)\mathsf{x}   \nonumber \\
	& =   e^{-  \gamma \nu}  \mathsf{x} \T (I-KC)\T  P^{-1}  \mathsf{x}     \nonumber \\
	& =  e^{-  \gamma \nu}  \mathsf{x}\T P^{-1} \mathsf{x}  - e^{-  \gamma \nu } \mathsf{x} \T C\T   (CPC\T +Q_t)^{-1} C   \mathsf{x}     \nonumber \\
	& \leq  e^{-\lambda^J_2} \mathcal{V}_2 , \quad \forall x_2\in \mathcal{J}_2   \label{eqn:V2+_2}
	\end{align}
	where  $\lambda^J_2: =   \gamma T_m  $, and we made use of the fact that $C\T   (CPC\T +Q_t)^{-1} C$ is positive semi-definite.

	Now, we are going to show the exponential stability of the set $ {\mathcal{A}}$ for overall hybrid closed-loop system $\mathcal{H}_2$ in (\ref{eqn:closed-loop2}).
	Let  $|x_2|_{\mathcal{A}} \geq 0$ be the distance of $x_2$ with respect to set $ {\mathcal{A}}$ such that $|x_2|_{ {\mathcal{A}}}^2:=\inf_{(\bar{R},\bar{\eta},\bar{\mathsf{x}},\bar{\tau})\in  {\mathcal{A}}}   (|\tilde{R}\bar{R}\T|_I^2 + \|\eta-\bar{\eta}\|^2 + \|\mathsf{x}- \bar{\mathsf{x}}\|^2 + \|\tau-\bar{\tau}\|^2 )=  \|\zeta\|^2 +   \|\mathsf{x}\|^2 $. Consider the Lyapunov function candidate $\mathcal{V}(x_2)= \varepsilon \mathcal{V}_1 (\tilde{R},\eta,\tau) + \mathcal{V}_2(\mathsf{x},\tau)$, with some $\varepsilon>0$. In view of (\ref{eqn:V1bound}) and (\ref{eqn:V2-bound_2}), one has
	\begin{align}
	& \underline{\alpha} |x_2|_{ {\mathcal{A}}}^2 \leq \mathcal{V}(x_2) \leq \bar{\alpha} |x_2|_{ {\mathcal{A}}}^2 \label{eqn:Vbound_2} 
	\end{align}
	where $\underline{\alpha}:=\min\{\varepsilon\underline{\alpha}_1, \underline{\alpha}_2'\}$ and $ \bar{\alpha}:=\max\{\varepsilon\bar{\alpha}_1,\bar{\alpha}_2'\}$. Applying the same steps as in (\ref{eqn:dotV}), from (\ref{eqn:V1})-(\ref{eqn:dotV1}) and  (\ref{eqn:V2-bound_2})-(\ref{eqn:dotV2_4})  one obtains
	\begin{equation}
	\dot{\mathcal{V}} = -\lambda_F \mathcal{V}, \quad \forall x_2\in \mathcal{F}_2  \label{eqn:dotV_2}
	\end{equation}
	with $	\varepsilon > \frac{\beta^2 p_M}{2\lambda_{1}^F \lambda_{2}^F \lambda_m^{P_1} e^{-  \gamma T_M}  }$ and $\lambda_F :=   \min\{ \frac{  \lambda_{1}^F  }{2} ,  ( \lambda_{2}^F -  \frac{\beta^2 p_M}{2\varepsilon \lambda_{1}^F \lambda_m^{P_1}  } )\}  $.  Moreover, in view of (\ref{eqn:V1+}) and (\ref{eqn:V2+_2}), one obtains
	\begin{equation}
	\mathcal{V}(x_2^+) 
	\leq e^{- \lambda_J} \mathcal{V}(x_2) , \quad \forall x_2\in \mathcal{J}_2 \label{eqn:V+_2}
	\end{equation}
	where $\lambda_J:= \min\{\lambda_1^J,\lambda_2^J\}$. Let $\lambda_c: = \min\{\lambda_F,\lambda_J\}$. In view of  (\ref{eqn:dotV_2}) and (\ref{eqn:V+_2}), one has $\mathcal{V}(x_2(t,j))\leq e^{-\lambda_c(t+j)} \mathcal{V}(x_2(0,0))$. Since $T_m$ is strictly positive by Assumption \ref{assum:1}, it follows that every maximal solution to the hybrid system $\mathcal{H}_2$ is complete. Then, applying (\ref{eqn:Vbound_2})  one can conclude that, for all $(t,j)\in \dom x_2$,
	$
	|x_2(t,j)|_{ {\mathcal{A}}} \leq \sqrt{\frac{ \bar{\alpha}}{\underline{\alpha}}}  e^{-\frac{1}{2}\lambda_c(t+j)}|x_2(0,0)|_{ {\mathcal{A}}},  
	$
	which shows that the set $ {\mathcal{A}}$ is exponentially stable. This completes the proof.


	\bibliographystyle{apalike}        
	\bibliography{autosam}           

\begin{thebibliography}{}

\bibitem[Alonge et~al., 2019]{alonge2019hybrid}
Alonge, F., D'Ippolito, F., Garraffa, G., and Sferlazza, A. (2019).
\newblock A hybrid observer for localization of mobile vehicles with
  asynchronous measurements.
\newblock {\em Asian Journal of Control}, 21(4):1506--1521.

\bibitem[Barrau and Bonnabel, 2017]{barrau2017invariant}
Barrau, A. and Bonnabel, S. (2017).
\newblock The invariant extended {Kalman} filter as a stable observer.
\newblock {\em IEEE Transactions on Automatic Control}, 62(4):1797--1812.

\bibitem[Berkane et~al., 2017]{berkane2017hybrid}
Berkane, S., Abdessameud, A., and Tayebi, A. (2017).
\newblock Hybrid attitude and gyro-bias observer design on {SO(3)}.
\newblock {\em IEEE Transactions on Automatic Control}, 62(11):6044--6050.

\bibitem[Berkane and Tayebi, 2017]{BerkaneCDC2017}
Berkane, S. and Tayebi, A. (2017).
\newblock Attitude observer using synchronous intermittent vector measurements.
\newblock In {\em Proc. 56th IEEE conference on decision and control}, pp.
  3027--3032.

\bibitem[Berkane and Tayebi, 2019]{berkane2019attitude}
Berkane, S. and Tayebi, A. (2019).
\newblock Attitude estimation with intermittent measurements.
\newblock {\em Automatica}, 105:415--421.

\bibitem[Burri et~al., 2016]{Burri25012016}
Burri, M., Nikolic, J., Gohl, P., Schneider, T., Rehder, J., Omari, S.,
  Achtelik, M., and Siegwart, R. (2016).
\newblock The {EuRoC} micro aerial vehicle datasets.
\newblock {\em The International Journal of Robotics Research},
  35(10):1157--1163.

\bibitem[Carnevale et~al., 2007]{carnevale2007lyapunov}
Carnevale, D., Teel, A.~R., and Nesic, D. (2007).
\newblock A lyapunov proof of an improved maximum allowable transfer interval
  for networked control systems.
\newblock {\em IEEE Transactions on Automatic Control}, 52(5):892--897.

\bibitem[Deyst and Price, 1968]{deyst1968conditions}
Deyst, J. and Price, C. (1968).
\newblock Conditions for asymptotic stability of the discrete minimum-variance
  linear estimator.
\newblock {\em IEEE Transactions on Automatic Control}, 13:702--705.

\bibitem[Ferrante et~al., 2016]{ferrante2016state}
Ferrante, F., Gouaisbaut, F., Sanfelice, R.~G., and Tarbouriech, S. (2016).
\newblock State estimation of linear systems in the presence of sporadic
  measurements.
\newblock {\em Automatica}, 73:101--109.

\bibitem[Goebel et~al., 2009]{goebel2009hybrid}
Goebel, R., Sanfelice, R., and Teel, A. (2009).
\newblock Hybrid dynamical systems.
\newblock {\em IEEE control systems magazine}, 29(2):28--93.

\bibitem[Goebel et~al., 2012]{goebel2012hybrid}
Goebel, R., Sanfelice, R., and Teel, A. (2012).
\newblock {\em Hybrid Dynamical Systems: modeling, stability, and robustness}.
\newblock Princeton University Press.

\bibitem[Grant et~al., 2009]{grant2009cvx}
Grant, M., Boyd, S., and Ye, Y. (2009).
\newblock Cvx: Matlab software for disciplined convex programming.

\bibitem[Hamel and Samson, 2018]{hamel2018riccati}
Hamel, T. and Samson, C. (2018).
\newblock Riccati observers for the nonstationary {PnP} problem.
\newblock {\em IEEE Transactions on Automatic Control}, 63(3):726--741.

\bibitem[Hesch et~al., 2013]{hesch2013consistency}
Hesch, J.~A., Kottas, D.~G., Bowman, S.~L., and Roumeliotis, S.~I. (2013).
\newblock Consistency analysis and improvement of vision-aided inertial
  navigation.
\newblock {\em IEEE Transactions on Robotics}, 30(1):158--176.

\bibitem[Hua and Allibert, 2018]{hua2018riccati}
Hua, M.-D. and Allibert, G. (2018).
\newblock Riccati observer design for pose, linear velocity and gravity
  direction estimation using landmark position and {IMU} measurements.
\newblock In {\em Proc. 2018 IEEE conference on Control Technology and
  Applications}, pp. 1313--1318. IEEE.

\bibitem[Hua et~al., 2015]{hua2015gradient}
Hua, M.-D., Hamel, T., Mahony, R., and Trumpf, J. (2015).
\newblock Gradient-like observer design on the {Special Euclidean group SE(3)}
  with system outputs on the real projective space.
\newblock In {\em Proc. 54th IEEE conference on decision and control}, pp.
  2139--2145.

\bibitem[Hua et~al., 2018]{hua2018attitude}
Hua, M.-D., Manerikar, N., Hamel, T., and Samson, C. (2018).
\newblock Attitude, linear velocity and depth estimation of a camera observing
  a planar target using continuous homography and inertial data.
\newblock In {\em Proc. IEEE International Conference on Robotics and
  Automation}, pp. 1429--1435. IEEE.

\bibitem[Jazwinski, 1970]{jazwinski1970}
Jazwinski, A.~H. (1970).
\newblock {\em Stochastic Processes and Filtering Theory}.
\newblock ACADEMIC PRESS, INC.

\bibitem[Kelly and Sukhatme, 2011]{kelly2011visual}
Kelly, J. and Sukhatme, G.~S. (2011).
\newblock Visual-inertial sensor fusion: Localization, mapping and
  sensor-to-sensor self-calibration.
\newblock {\em The International Journal of Robotics Research}, 30(1):56--79.

\bibitem[Khosravian et~al., 2015]{khosravian2015observers}
Khosravian, A., Trumpf, J., Mahony, R., and Lageman, C. (2015).
\newblock Observers for invariant systems on {Lie} groups with biased input
  measurements and homogeneous outputs.
\newblock {\em Automatica}, 55:19--26.

\bibitem[Li et~al., 2017]{li2017robust}
Li, Y., Phillips, S., and Sanfelice, R.~G. (2017).
\newblock Robust distributed estimation for linear systems under intermittent
  information.
\newblock {\em IEEE Transactions on Automatic Control}, 63(4):973--988.

\bibitem[Mourikis and Roumeliotis, 2007]{mourikis2007multi}
Mourikis, A.~I. and Roumeliotis, S.~I. (2007).
\newblock A multi-state constraint {Kalman} filter for vision-aided inertial
  navigation.
\newblock In {\em Proc. of IEEE International Conference on Robotics and
  Automation (ICRA)}, pp. 3565--3572.

\bibitem[Mourikis et~al., 2009]{mourikis2009vision}
Mourikis, A.~I., Trawny, N., Roumeliotis, S.~I., Johnson, A.~E., Ansar, A., and
  Matthies, L. (2009).
\newblock Vision-aided inertial navigation for spacecraft entry, descent, and
  landing.
\newblock {\em IEEE Transactions on Robotics}, 25(2):264--280.

\bibitem[Rehbinder and Ghosh, 2003]{rehbinder2003pose}
Rehbinder, H. and Ghosh, B.~K. (2003).
\newblock Pose estimation using line-based dynamic vision and inertial sensors.
\newblock {\em IEEE Transactions on Automatic Control}, 48(2):186--199.

\bibitem[Scaramuzza and Fraundorfer, 2011]{scaramuzza2011visual}
Scaramuzza, D. and Fraundorfer, F. (2011).
\newblock Visual odometry [tutorial].
\newblock {\em IEEE robotics \& automation magazine}, 18(4):80--92.

\bibitem[Sferlazza et~al., 2019]{sferlazza2019time}
Sferlazza, A., Tarbouriech, S., and Zaccarian, L. (2019).
\newblock Time-varying sampled-data observer with asynchronous measurements.
\newblock {\em IEEE Transactions on Automatic Control}, 64(2):869--876.

\bibitem[Shi and Tomasi, 1994]{shi1994good}
Shi, J. and Tomasi, C. (1994).
\newblock Good features to track.
\newblock In {\em Proc. IEEE conference on Computer Vision and Pattern
  Recognition}, pp. 593--600. IEEE.

\bibitem[Tayebi et~al., 2013]{tayebi2013inertial}
Tayebi, A., Roberts, A., and Benallegue, A. (2013).
\newblock Inertial vector measurements based velocity-free attitude
  stabilization.
\newblock {\em IEEE Transactions on Automatic Control}, 58(11):2893--2898.

\bibitem[Teel et~al., 2013]{teel2013lyapunov}
Teel, A.~R., Forni, F., and Zaccarian, L. (2013).
\newblock Lyapunov-based sufficient conditions for exponential stability in
  hybrid systems.
\newblock {\em IEEE Transactions on Automatic Control}, 58(6):1591--1596.

\bibitem[Vasconcelos et~al., 2010]{vasconcelos2010nonlinear}
Vasconcelos, J., Cunha, R., Silvestre, C., and Oliveira, P. (2010).
\newblock A nonlinear position and attitude observer on {SE}(3) using landmark
  measurements.
\newblock {\em Systems \& Control Letters}, 59(3-4):155--166.

\bibitem[Wang and Tayebi, 2017]{wang2017globally}
Wang, M. and Tayebi, A. (2017).
\newblock Globally asymptotically stable hybrid observers design on {SE}(3).
\newblock In {\em Proc. 56th IEEE conference on decision and control}, pp.
  3033--3038.

\bibitem[Wang and Tayebi, 2018]{wang2018navigation}
Wang, M. and Tayebi, A. (2018).
\newblock A globally exponentially stable nonlinear hybrid observer for {3D}
  inertial navigation.
\newblock In {\em Proc. 57th IEEE conference on decision and control}, pp.
  1367--1372.

\bibitem[Wang and Tayebi, 2019]{wang2019hybrid}
Wang, M. and Tayebi, A. (2019).
\newblock Hybrid pose and velocity-bias estimation on {SE(3)} using inertial
  and landmark measurements.
\newblock {\em IEEE Transactions on Automatic Control}, 64(8):3399--3406.

\bibitem[Wang and Tayebi, 2020]{wang2020hybrid}
Wang, M. and Tayebi, A. (2020).
\newblock Hybrid nonlinear observers for inertial navigation using landmark
  measurements.
\newblock {\em IEEE Transactions on Automatic Control}.
\newblock doi:
  \href{https://doi.org/10.1109/TAC.2020.2972213}{10.1109/TAC.2020.2972213}.

\end{thebibliography}

\end{document}